\documentstyle[twoside,11pt]{article}

\begin{document}
\title{{\bf Convergence Rate of Zero Viscosity Limit on Large Amplitude Solution
to a Conservation Laws Arising in Chemotaxis}}
\author{H{\sc ongyun} P{\sc eng}, L{\sc izhi} R{\sc uan}\thanks{Corresponding author.\ \ Email:
rlz@mail.ccnu.edu.cn}, C{\sc hangjiang} Z{\sc hu}\\
The Hubei Key Laboratory of Mathematical Physics\\
School of Mathematics and Statistics\\
Central China Normal University, Wuhan, 430079, P. R. China}

\date{}
\maketitle

\tableofcontents

\arraycolsep=1.5pt
\newtheorem{Lemma}{Lemma}[section]
\newtheorem{Theorem}{Theorem}[section]
\newtheorem{Definition}{Definition}[section]
\newtheorem{Proposition}{Proposition}[section]
\newtheorem{Remark}{Remark}[section]
\newtheorem{Corollary}{Corollary}[section]

\begin{abstract}
In this paper, we investigate large amplitude solutions to a system
of conservation laws which is transformed, by a change of variable,
from the well-known Keller-Segel model describing cell (bacteria)
movement toward the concentration gradient of the chemical that is
consumed by the cells. For the Cauchy problem and initial-boundary
value problem, the global unique solvability is proved based on the
energy method. In particular, our main purpose is to investigate the
convergence rates as the diffusion parameter $\varepsilon$ goes to
zero. It is shown that the convergence rates in $L^\infty$-norm are
of the order $O\left(\varepsilon\right)$ and $O(\varepsilon^{3/4})$
corresponding to the Cauchy problem and the initial-boundary value
problem respectively.

\vspace{4mm}



\end{abstract}

\section{Introduction}
In this paper, we investigate large amplitude solutions to a system
of conservation laws
\begin{equation}\label{Equation}
\left\{
\begin{array}{l}
u_t^\varepsilon+\left(\varepsilon \left(u^\varepsilon\right)^2-v^\varepsilon\right)_x
=\varepsilon u_{xx}^\varepsilon,\\[3mm]
v_t^\varepsilon-\left(u^\varepsilon
v^\varepsilon\right)_x=v_{xx}^\varepsilon,
\end{array}
\right.
\end{equation}
which is transformed, by a change of variable, from the well-known
Keller-Segel model describing cell (bacteria) movement toward the
concentration gradient of the chemical that is consumed by the
cells. Here $\varepsilon>0$ is a positive constant.

The aim here is to study the global unique solvability on the Cauchy
problem and initial-boundary value problem in the framework of
large-amplitude $H^2$ solutions. In particular, our main purpose is
to investigate the convergence rates as the diffusion parameter
$\varepsilon$ goes to zero. It is shown that the convergence rates
in $L^\infty$-norm are of the order $O\left(\varepsilon\right)$ and
$O(\varepsilon^{3/4})$ corresponding to the Cauchy problem and the
initial-boundary value problem respectively.

Firstly we are concerned with the Cauchy problem of (\ref{Equation})
with initial data
\begin{equation}\label{Initial}
\left(u^\varepsilon,
v^\varepsilon\right)(x,0)=\left(u_0^\varepsilon(x),v_0^\varepsilon(x)\right)\rightarrow(0,v_\infty)\
\ {\rm as}\ \ x\rightarrow\pm\infty.
\end{equation}
Hereafter, $v_\infty$ is a given positive constant. On one hand, the
global existence of solutions to the Cauchy problem
(\ref{Equation}), (\ref{Initial}) for any fixed $\varepsilon>0$ is
shown. On the other hand, we prove that solutions to the Cauchy
problem (\ref{Equation}), (\ref{Initial}) convergence to solutions
to the problem of zero viscosity limit when $\varepsilon\rightarrow
0^+$, which is formally formulated as follows:
\begin{equation}\label{Limit equation}
\left\{
\begin{array}{l}
u_t^0-v_x^0=0,\\[3mm]
v_t^0-\left(u^0 v^0\right)_x=v_{xx}^0,
\end{array}
\right.
\end{equation}
with initial data
\begin{equation}\label{Limit initial}
\left(u^0,
v^0\right)(x,0)=\left(u_0^0(x),v_0^0(x)\right)\rightarrow(0,v_\infty)\
\ {\rm as}\ \ x\rightarrow\pm\infty.
\end{equation}

In addition, we consider the initial-boundary value problem of
(\ref{Equation}) with initial data
\begin{equation}\label{Initial1}
\left(u^\varepsilon,
v^\varepsilon\right)(x,0)=\left(u_0^\varepsilon(x),v_0^\varepsilon(x)\right),\
\ \ 0\leq x \leq 1,
\end{equation}
and the boundary conditions
\begin{equation}\label{boundary}
\left(u^\varepsilon,
v^\varepsilon_x\right)(0,t)=\left(u^\varepsilon,v_x^\varepsilon\right)(1,t)=(0,0),\
\ \ t\geq 0,
\end{equation}
which implies
\begin{equation}\label{1boundary}
u_{xx}^\varepsilon(0,t)=u_{xx}^\varepsilon(1,t)=0,\ \ \ t\geq 0,
\end{equation}
where the compatibility conditions
$u_0^\varepsilon(0)=u_0^\varepsilon(1)=0$ and
$v_{0x}^\varepsilon(0)=v_{0x}^\varepsilon(1)=0$ are satisfied.

{} For one thing, the global existence of solutions to the
initial-boundary value problem of (\ref{Equation}),
(\ref{Initial1}), (\ref{boundary}) for any fixed $\varepsilon>0$ is
shown. For another thing, we prove that solutions to the
initial-boundary value problem of (\ref{Equation}),
(\ref{Initial1}), (\ref{boundary}) convergence to solutions to the
problem of zero viscosity limit when $\varepsilon\rightarrow 0^+$,
which is formally formulated as follows:
\begin{equation}\label{Limit equation1}
\left\{
\begin{array}{l}
u_t^0-v_x^0=0,\\[3mm]
v_t^0-\left(u^0 v^0\right)_x=v_{xx}^0,
\end{array}
\right.
\end{equation}
with initial data
\begin{equation}\label{Initial2}
\left(u^0, v^0\right)(x,0)=\left(u_0^0(x),v_0^0(x)\right),\ \ \
0\leq x \leq 1,
\end{equation}
and the boundary conditions
\begin{equation}\label{boundary1}
u^0(0,t)=u^0(1,t)=0, \ \ t\geq0,
\end{equation}
which implies
\begin{equation}\label{boundary11}
v_x^0(0,t)=v_x^0(1,t)=0, \ \ t\geq0,
\end{equation}
where the compatibility conditions $u_0^0(0)=u_0^0(1)=0$ and
$v_{0x}^0(0)=v_{0x}^0(1)=0$ are satisfied.

Precisely speaking, these results are given as follows. We have the
first main theorem, which is concerned with the global existence of
the Cauchy problem (\ref{Equation})-(\ref{Initial}).

\begin{Theorem}[Global existence]
Let $\left(u_0^\varepsilon,v_0^\varepsilon-v_\infty\right)\in
H^2({\bf R})$ and assume that there exists a positive constant
$\alpha>0$ such that $\inf\limits_{x\in{\bf
R}}v_0^\varepsilon(x)\geq \alpha$. Then there exists a unique global
solution $\left(u^\varepsilon(x,t),v^\varepsilon(x,t)\right)$ of the
Cauchy problem (\ref{Equation})-(\ref{Initial}) which satisfies
\begin{equation}\label{Existence1-c}
\left\{
\begin{array}{l}
\left(u^\varepsilon, v^\varepsilon-v_\infty\right) \in
L^\infty([0,\infty),H^2({\bf R})),\\[3mm]
 \frac{v_x^\varepsilon}{\sqrt{v^\varepsilon}}\in L^2([0,\infty),L^2({\bf R})),\quad
 v_{xx}^\varepsilon\in L^2([0,\infty),H^1({\bf R}))
 \end{array}
 \right.
\end{equation}
and
\begin{equation}\label{Existence3-c}
\begin{array}{ll}
\displaystyle\left\|u^\varepsilon(t)\right\|^2_{H^2}+\left\|v^\varepsilon(t)-v_\infty\right\|^2_{H^2}+\varepsilon\int_0^t
\|u_x^\varepsilon(s)\|^2_{H^2}ds\\[3mm]
\displaystyle\hspace{1.8cm}+\int_0^t
\left\{\left\|\frac{v_x^\varepsilon(s)}{\sqrt{v^\varepsilon}}\right\|^2_{L^2}
+\left\|\sqrt{v^\varepsilon}u_x^\varepsilon(s)\right\|^2_{L^2}
+\|v_{xx}^\varepsilon(s)\|^2_{H^1}\right\}ds\\[3mm]
\displaystyle\hspace{6.5cm}\leq
C\left(\left\|u_0^\varepsilon\right\|^2_{H^2}+\left\|v_0^\varepsilon-v_\infty\right\|^2_{H^2}\right).
\end{array}
\end{equation}
Moreover, for any fixed $T>0$, there exists a positive constant
$C(T)$ which depends only on $T$ such that
\begin{eqnarray}\label{Existence2-c}
v^\varepsilon(x,t)\geq C(T).
\end{eqnarray}
Here $C$ and $C(T)$ are all positive constants independent of
$\varepsilon$.
\end{Theorem}

\noindent\textbf{Remark 1.1.}\ \ In the case of $\varepsilon=0$, the
global unique solvability on the Cauchy problem (\ref{Limit
equation})-(\ref{Limit initial}) was included in \cite{Guo09}, which
will be summarized in Lemma 2.4 later. As is well-known that it is
more difficult when the nonlinearity is higher in the setting of
large amplitude. It is nontrivial to obtain the estimate on
(\ref{Existence3-c}) in the setting of large amplitude for the case
of $\varepsilon>0$ since the nonlinear term
$\left\{\left(\varepsilon u^\varepsilon\right)^2\right\}_x$ appears.
We overcome these difficulties by applying some suitable
interpolation inequalities.

Next, the following main result refers to convergence rate under the
additional assumptions on initial data.

\begin{Theorem}[$L^\infty$-convergence rate]
Under the same conditions of Theorem 1.1, we assume the initial data
$\left(u_0^\varepsilon,v_0^\varepsilon\right)$ and
$\left(u_0^0,v_0^0\right)$ satisfy
$\left\|u_0^\varepsilon-u_0^0\right\|_{H^2({\bf
R})}+\left\|v_0^\varepsilon-v_0^0\right\|_{H^2({\bf R})}
=O(\varepsilon)$, and the further regularity
$\left(u_0^0,v_0^0\right)\in H^3({\bf R})$. Then we have
\begin{eqnarray}\label{Th-con-c}
\arraycolsep=1.5pt
\begin{array}[b]{rl}\displaystyle
\left\|\left(u^\varepsilon-u^0\right)(t)\right\|_{L^\infty({\bf R})}
+\left\|\left(v^\varepsilon-v^0\right)(t)\right\|_{L^\infty({\bf
R})}\leq&\displaystyle C\varepsilon.
\end{array}
\end{eqnarray}
Here $C>0$ is a positive constant independent of $\varepsilon$.
\end{Theorem}

\noindent\textbf{Remark 1.2.}\ \ In the case of $\varepsilon>0$, for
the initial-boundary value problem of (\ref{Equation}) with
Dirichlet boundary conditions
$u^\varepsilon(0,t)=u^\varepsilon(1,t)=0,\ \
v^\varepsilon(0,t)=v^\varepsilon(1,t)=0$, with the help of the
methods in \cite{Zhang07}, we can prove the global unique
solvability in the setting of small amplitude since the difficulties
derived from the boundary effect on the viscosity $\varepsilon
u_{xx}^\varepsilon$. What's more, the boundary layer will occur and
the boundary layer thickness will be shown as in
\cite{Frid99,Jiang09,Ruan11}. These are left to our future study
since these are not in accordance with our main aim in the paper,
which mainly refers to large amplitude solution.

Our final result focusing on the initial-boundary value problem is
stated as follows.
\begin{Theorem}[Initial-boundary value problem]
Let $\left(u_0^\varepsilon,v_0^\varepsilon\right)\in H^2({[0,1]})$
and assume that there exists a positive constant $\alpha>0$ such
that $\inf\limits_{x\in{[0,1]}}v_0^\varepsilon(x)\geq \alpha$.
\end{Theorem}
({\rm i})\ \ Then there exists a unique global solution
$\left(u^\varepsilon(x,t),v^\varepsilon(x,t)\right)$ of the
initial-boundary value problem (\ref{Equation}), (\ref{Initial1})
and (\ref{boundary}) which satisfies
\begin{equation}\label{Existence1-cb}
\left\{
\begin{array}{l}
\left(u^\varepsilon, v^\varepsilon\right) \in
L^\infty([0,\infty),H^2([0,1])),\\[3mm]
 \frac{v_x^\varepsilon}{\sqrt{v^\varepsilon}}\in L^2([0,\infty),L^2([0,1])),\quad
 v_{xx}^\varepsilon\in L^2([0,\infty),H^1([0,1]))
 \end{array}
 \right.
\end{equation}
and
\begin{equation}\label{Existence3-cb}
\begin{array}{ll}
\displaystyle\left\|u^\varepsilon(t)\right\|^2_{H^2}+\left\|v^\varepsilon(t)\right\|^2_{H^2}+\varepsilon\int_0^t
\|u_x^\varepsilon(s)\|^2_{H^2}ds\\[3mm]
\displaystyle\hspace{1.8cm}+\int_0^t
\left\{\left\|\frac{v_x^\varepsilon(s)}{\sqrt{v^\varepsilon}}\right\|^2_{L^2}
+\left\|\sqrt{v^\varepsilon}u_x^\varepsilon(s)\right\|^2_{L^2}
+\|v_{xx}^\varepsilon(s)\|^2_{H^1}\right\}ds\\[3mm]
\displaystyle\hspace{6.5cm}\leq
C\left(\left\|u_0^\varepsilon\right\|^2_{H^2}+\left\|v_0^\varepsilon\right\|^2_{H^2}\right).
\end{array}
\end{equation}
Moreover, for any fixed $T>0$, there exists a positive constant
$C(T)$ which depends only on $T$ such that
\begin{eqnarray}\label{Existence2-cb}
v^\varepsilon(x,t)\geq C(T).
\end{eqnarray}
({\rm ii})\ \ We assume the initial data
$\left(u_0^\varepsilon,v_0^\varepsilon\right)$ and
$\left(u_0^0,v_0^0\right)$ satisfy
$\left\|u_0^\varepsilon-u_0^0\right\|_{H^2([0,1])}+\left\|v_0^\varepsilon-v_0^0\right\|_{H^2([0,1])}
=O(\varepsilon)$. Then we have
\begin{eqnarray}\label{Th-con-cb}
\arraycolsep=1.5pt
\begin{array}[b]{rl}\displaystyle
\left\|\left(u^\varepsilon-u^0\right)(t)\right\|_{L^\infty([0,1])}
+\left\|\left(v^\varepsilon-v^0\right)(t)\right\|_{L^\infty([0,1])}\leq&\displaystyle
C\varepsilon^{\frac{3}{4}}.
\end{array}
\end{eqnarray}
Here $C$ and $C(T)$ are all positive constants independent of
$\varepsilon$.

\noindent\textbf{Remark 1.3.}\ \ Convergence rates (\ref{Th-con-cb})
on the initial-boundary value problem is slower than one in
(\ref{Th-con-c}) on the Cauchy problem, which is caused by the
boundary effect. Essentially it is different between
initial-boundary value problem and Cauchy problem that the
regularity of solutions to the Cauchy problem can be improved when
the regularity of initial data is imposed further regularity whereas
it fails for the initial-boundary value problem.

The problem of the zero viscosity limit is one of the important
topics. In particular, when parabolic equations with small viscosity
are applied as perturbations, convergence rate of the Cauchy problem
or the boundary layer question of the initial-boundary value problem
for many other equations also arises in the theory of hyperbolic
systems in the case of one-dimension or multi-dimension, cf.
\cite{Chen08,Frid99,Jiang09,Ruan11,Xin98,Xin99}. To our knowledge,
fewer results on the equations (\ref{Equation}) have been obtained
in this direction.

Now let us review some known results related to the system
(\ref{Equation}) and (\ref{Limit equation}), which has been
extensively studied by several authors in different contexts, cf.
\cite{Guo09,Li111,Li11,Zhang07}. The conservation laws
(\ref{Equation}) are derived from the original well-known
Keller-Segel model
\begin{equation}\label{Keller-Segel}
\left\{
\begin{array}{l}
c_t=\varepsilon c_{xx}-uf(c),\\[3mm]
u_t=\left(Du_x-\chi uc^{-1}c_x\right)_x,
\end{array}
\right.
\end{equation}
which was proposed by Keller and Segel in \cite{Keller71} to
describe the traveling band behavior of bacteria due to the
chemotactic response (i.e., the oriented movement of cells to the
chemical concentration gradient) observed in experiments
\cite{Adler66,Adler69}. In model (\ref{Keller-Segel}), $u(x,t)$ and
$c(x,t)$ denote the cell density and the chemical concentration,
respectively. $D>0$ is the diffusion rate of cells (bacteria) and
$\varepsilon>0$ is the diffusion rate of
 chemical substance. $\chi$ is a positive constant often referred
to as chemosensitivity. $f(c)$ is a kinetic function describing the
chemical reaction between cells and the chemical.

When $f(c)$ is a positive constant, namely, $f(c)=\alpha>0$, the
existence of traveling wave solutions of (\ref{Keller-Segel}) with
$\varepsilon=0$ was established by Keller and Segel themselves in
\cite{Keller71}. When $\varepsilon\neq 0$, the existence and linear
instability of traveling wave solutions of (\ref{Keller-Segel}) were
shown by Nagai and Ikeda in \cite{Nagai91} where the authors also
obtained the diffusion limits of traveling wave solutions of
(\ref{Keller-Segel}) as $\varepsilon$ approaches zero. Precisely
they proved that the traveling wave solution in the form
$\left(B^\varepsilon(x,t),
S^\varepsilon(x,t)\right)=\left(B^\varepsilon(z),
S^\varepsilon(z)\right) (z=x-ct)$ approximates to the corresponding
to traveling wave solution $\left(B^0(x,t),
S^0(x,t)\right)=\left(B^0(z), S^0(z)\right) (z=x-ct)$ when
$\varepsilon$ goes to zero.  For the reduction of Keller-Segel
system (\ref{Keller-Segel}) to system (\ref{Equation}), we refer to
Refs. \cite{Li11}, see also Appendix in this paper.

In \cite{Li11}, Li and Wang established the existence and the
nonlinear stability of traveling wave solutions to a system of
conservation laws (\ref{Equation}). They prove the existence of
traveling fronts by the phase plane analysis and show the asymptotic
nonlinear stability of traveling wave solutions without the
smallness assumption on the wave strengths by the method of energy
estimates.

  There is limiting case of the Keller-Segel model (\ref{Keller-Segel}).
That is when the diffusion of chemical substance is so small that it
is negligible, i.e, $\varepsilon\rightarrow0^+$, then the model
(\ref{Keller-Segel}) becomes

\begin{equation}\label{Othmer-Stevens}
\left\{
\begin{array}{l}
c_t=-uf(c),\\[3mm]
u_t=\left(Du_x-\chi uc^{-1}c_x\right)_x.
\end{array}
\right.
\end{equation}

A version of system (\ref{Othmer-Stevens}) was proposed by Othmer
and Stevens in \cite{Othmer97} to describe the chemotactic movement
of particles where the chemicals are non-diffusible. Othmer and
Stevens in \cite{Othmer97} have developed a number of mathematical
models of chemotaxis to illustrate aggregation leading (numerically)
to nonconstant steady-states, blow-up resulting in the formation of
singularities and collapse or the formation of a spatially uniform
steady state. The models developed in \cite{Othmer97} have been
studied in depth by Levine and Sleeman in \cite{Levine97}. They gave
some heuristic understanding of some of these phenomena and
investigated the properties of solutions of a system of chemotaxis
equation arising in the theory of reinforced random walks. Y. Yang,
H. Chen and W.A. Liu in \cite{Yang01} studied the global existence
and blow-up in a finite-time of solutions for the case considered in
\cite{Levine97}, respectively. They found that even at the same
growth rate the behavior of the biological systems can be very
different just because they started their action in different
conditions. For the other results on the initial-boundary value
problem of (\ref{Othmer-Stevens}) refer to \cite{Hillen04,
Levine01}.

Recently, the modified model related closely to chemotaxis was also
investigated by Painter and Hillen in \cite{Painter11}. They
explored the dynamics of a one-dimensional Keller-Segel type model
for chemotaxis incorporating a logistic cell growth term.

Finally, we have to mention the work in \cite{Zhang07, Guo09}, which
investigated the initial-boundary value problem and Cauchy problem
on limit equation (\ref{Limit equation}) and motivate our
investigation in this paper. Zhang and Zhu in \cite{Zhang07} studied
the initial-boundary value problem of (\ref{Limit equation}) with
initial data
\begin{equation}\label{Initial-Zhang}
\left(u^0, v^0\right)(x,0)=\left(u_0^0(x),v_0^0(x)\right)
\end{equation}
and boundary conditions
\begin{equation}\label{Boundary-Zhang}
u^0(0,t)=u^0(1,t)=0.
\end{equation}
Let $\left(u_0^0,v_0^0-1\right)\in H^2([0,1])$ and assume that there
exists a positive constant $\epsilon>0$ sufficiently small such that
$\left\|u_0^0\right\|^2_{H^2}+\left\|v_0^0-1\right\|^2_{H^2}\leq\epsilon$.
Then there exists a unique global solution
$\left(u^0(x,t),v^0(x,t)\right)$ of the initial-boundary value
problem (\ref{Limit equation}), (\ref{Initial-Zhang}) and
(\ref{Boundary-Zhang}) which satisfies
\begin{eqnarray}\label{existence1-Zhang}
\left(u^0, v^0\right) \in L^\infty([0,\infty),H^2([0,1])),\quad
 v_x^0\in L^2([0,\infty),H^2([0,1])).
\end{eqnarray}

Guo etc. in \cite{Guo09} studied the Cauchy problem of (\ref{Limit
equation}) with initial data
\begin{equation}\label{Initial-Guo}
\left(u^0,
v^0\right)(x,0)=\left(u_0^0(x),v_0^0(x)\right)\rightarrow(u_\pm,
v_\infty)\ \ {\rm as}\ \ x\rightarrow\pm\infty.
\end{equation}
They proved the existence of global solutions to the Cauchy problem
of a hyperbolic-parabolic coupled system with large initial data,
which is summarized in Lemma 2.4 later.

\bigbreak\noindent\textbf{Notations:}\ \ Throughout this paper, we
denote positive constants by $C$ may depends on $T$, but is
independent of $\varepsilon$. Moreover, the character ``$C$" may
differ in different places. $L^p=L^p({\Omega})$ $(1\leq
p\leq\infty)$ denotes usual Lebesgue space with the norm
$$
\begin{array}{c}
\bigbreak\displaystyle \|f\|_{L^p({\Omega})}=\left(\int_{\Omega}
|f(x)|^pdx\right)^{\frac{1}{p}}, \ \
1\leq p<\infty,\\
\|f\|_{L^\infty}=\sup\limits_{x\in {\Omega}}|f(x)|.
\end{array}
$$
$H^l({\Omega}) \ (l\geq 0)$ denotes the usual $l$th-order Sobolev
space with the norm
$$
\|f\|_l=\left(\sum\limits^l_{j=0}\|\partial^j_xf\|^2\right)^{\frac{1}{2}},
$$
where $\Omega={\bf R}$ or $[0,1]$, and
$\|\cdot\|=\|\cdot\|_0=||\cdot\|_{L^2}$. For simplicity, $\|f(\cdot,
t)\|_{L^p}$ and $\|f(\cdot, t)\|_l$ are denoted by $\|f(t)\|_{L^p}$
and $\|f(t)\|_{l}$ respectively.

The rest of the paper is organized as follows. In Section 2, we
study the Cauchy problem (\ref{Equation}) and (\ref{Initial}). For
any fixed $\varepsilon>0$, the global unique solvability are
obtained. When diffusion parameter $\varepsilon\rightarrow0^+$, the
limit problem is considered and we show that the convergence rates
in $L^\infty$-norm is of the order $O(\varepsilon)$. In Section 3,
we study the initial-boundary value problem (\ref{Equation}),
(\ref{Initial1}) and (\ref{boundary}).  In the last Section 4, for
the convenience of the readers, we use the appendix to give the
derivation of system (\ref{Equation}) and (\ref{Limit equation}).

\section{Cauchy problem}
\subsection{Global existence on the case of $\varepsilon>0$}
\setcounter{equation}{0} In this section, we are concerned with the
global existence of large-amplitude $H^2$ solutions to the Cauchy
problem (\ref{Equation}), (\ref{Initial}), which will be proven by
continuing a unique local solution using a uniform-in-$\varepsilon$
\textit{a priori} estimate in the setting of large data. The
construction on the local existence of the solutions is standard
based on iteration argument and Fixed Point Theorem. Next, we will
devoted ourself to obtaining some uniform-in-$\varepsilon$ \textit{a
priori} estimates on the solution $\left(u^\varepsilon(x,t),
v^\varepsilon(x,t)\right)$ to the Cauchy problem (\ref{Equation}),
(\ref{Initial}) for all $\varepsilon>0$. To this end, $
v^\varepsilon(x,t)$ is supposed to satisfy the priori assumption
\begin{equation}\label{priori assumption-c}
v^\varepsilon(x,t)>0,\ \ (x,t)\in{\bf R}\times[0,T]
\end{equation}
{}for any fixed $T>0$.

\begin{Lemma} [Basic energy estimate]
Assume that the assumptions listed in Theorem 1.1 are satisfied.
Then there exists a positive constant $C$, independent of
$\varepsilon$, such that
\begin{eqnarray}\label{Le-l.basic1-c}
\int_{\bf R} \eta dx+\int_0^t\int_{\bf
R}{\frac{\left(v_x^\varepsilon\right)^2 } {v^\varepsilon}}dxds
+\varepsilon\int_0^t\int_{\bf R} \left(u_x^\varepsilon\right)^2dxds=
\int_{\bf R}{\eta _0 } dx.
\end{eqnarray}
\end{Lemma}

\noindent\textbf{Proof.}\ \ As in \cite{Zhang07, Guo09}, it is easy
to see that
\begin{eqnarray*}
\left\{
\begin{array}{rl}
  \eta\left(u^\varepsilon,v^\varepsilon\right) =& \displaystyle\frac{1} {2}\left(u^\varepsilon\right)^2  +
  v^\varepsilon\ln \frac{v^\varepsilon}
{{v_\infty}} - \left(v^\varepsilon - v_\infty\right), \\[3mm]
  q\left(u^\varepsilon,v^\varepsilon\right) =& \displaystyle - u^\varepsilon v^\varepsilon\ln \frac{v^\varepsilon}
{{v_\infty}}+\frac{2}{3}\varepsilon \left(u^\varepsilon\right)^3
\end{array}
\right.
\end{eqnarray*}
is an entropy-entropy flux pair to the hyperbolic system
(\ref{Equation}) normalized at $(0, v_\infty)$. Moreover, based on
the definition of the entropy-entropy flux pair (see
\cite{Smoller83}), it is easy to verify that
$\left(\eta\left(u^\varepsilon,v^\varepsilon\right),q(u^\varepsilon,v^\varepsilon)\right)$
satisfies
\begin{equation}\label{Le-l.basicp1-c}
\eta _t  + q_x  =\varepsilon u^\varepsilon
u_{xx}^\varepsilon+v_{xx}^\varepsilon \ln \frac{v^\varepsilon}
{{v_\infty}}.
\end{equation}
Integrating the above identity with respect $x$ and $t$ over ${\bf
R}\times [0,t]$, by using integration by parts we can get
(\ref{Le-l.basic1-c}). This completes the proof of Lemma 2.1.


\begin{Lemma} [First order energy estimate]
Assume that the assumptions listed in Theorem 1.1 are satisfied.
Then there exists a positive constant $C$, independent of
$\varepsilon$, such that
\begin{eqnarray}\label{Le-v-first1-c}
\begin{array}{ll}
\left\| u^\varepsilon\right\|_{L^\infty ({\bf R}\times [0,T])
}&\displaystyle+\left\| v^\varepsilon \right\|_{L^\infty ({\bf
R}\times [0,T])
}+\int_{\bf R}\left\{\left(u_x^\varepsilon\right)^2+\left(v_x^\varepsilon\right)^2\right\}dx\\[3mm]
&\displaystyle+\int^t_0\int_{\bf R} \left\{v^\varepsilon
\left(u_x^\varepsilon\right)^2+\left(v_{xx}^\varepsilon\right)^2\right\}dxds+\varepsilon\int^t_0\int_{\bf
R}\left(u_{xx}^\varepsilon\right)^2dxds\leq C.
\end{array}
\end{eqnarray}
\end{Lemma}
Here and hereafter, $C$ denotes positive constant depends on $T$.

 \noindent\textbf{Proof.}\ \ As in \cite{Guo09},  notice that
\begin{equation}\label{Le-v-firstp1-c}
u_{xt}^\varepsilon =
v_{xx}^\varepsilon-\varepsilon\left[\left(u^\varepsilon\right)^2\right]_{xx}+\varepsilon
u_{xxx}^\varepsilon= v_t^\varepsilon-\left(u^\varepsilon
v^\varepsilon\right)_x-\left[\varepsilon\left(u^\varepsilon\right)^2
\right]_{xx}+\varepsilon u_{xxx}^\varepsilon.
\end{equation}
Multiplying (\ref{Le-v-firstp1-c}) by $2u_x^\varepsilon$, then
integrating the resulting equation with respect $x$ and $t$ over
${\bf R}\times {[0,t]}$, we have by exploiting some integrations by
parts and Cauchy-Schwarz inequality that
\begin{eqnarray}\label{Le-v-firstp2-c}
\arraycolsep=1.5pt
\begin{array}[b]{rl}\displaystyle
 &\displaystyle\int_{\bf R}\left(u_x^\varepsilon\right)^2dx
 - \int_{\bf R}\left(u_{0x}^\varepsilon\right)^2dx+2\int_0^t\int_{\bf R}v^\varepsilon \left(u_x^\varepsilon\right)^2dxds
 +2\varepsilon\int_0^t\int_{\bf R}\left(u_{xx}^\varepsilon\right)^2dxds\\[3mm]
 =&\displaystyle 2\int_0^t\int_{\bf R}v_t^\varepsilon u_x^\varepsilon dxds
  - 2\int_0^t\int_{\bf R}u^\varepsilon v_x^\varepsilon u_x^\varepsilon dxds
  -2\varepsilon\int_0^t\int_{\bf R}\left(u_x^\varepsilon\right)^3dxds \\[3mm]
  = &\displaystyle \sum\limits_{i=1}^3I_i.
 \end{array}
\end{eqnarray}
$I_1$-$I_3$ are estimated as follows:
\begin{eqnarray}\label{Le-v-firstp2-I1-c}
\arraycolsep=1.5pt
\begin{array}[b]{rl}\displaystyle
 I_1=&\displaystyle-2\int_{\bf R} \left(v_0^\varepsilon- v_\infty\right)u_{0x}^\varepsilon dx+
 2\int_{\bf R} \left(v^\varepsilon- v_\infty\right)u_x^\varepsilon dx
  - 2\int_0^t\int_{\bf R}\left(v^\varepsilon- v_\infty\right)u_{xt}^\varepsilon dxds\\[3mm]
 =  &\displaystyle-2\int_{\bf R} \left(v_0^\varepsilon- v_\infty\right)u_{0x}^\varepsilon dx
 +2\int_{\bf R} \left(v^\varepsilon- v_\infty\right)u_x^\varepsilon dx
  - 2\int_0^t\int_{\bf R}\left(v^\varepsilon- v_\infty\right)v_{xx}^\varepsilon dxds\\[3mm]
 &\displaystyle+2\varepsilon\int_0^t\int_{\bf R}\left(v^\varepsilon- v_\infty\right)
 \left[\left(u^\varepsilon\right)^2\right]_{xx} dxds
 - 2\varepsilon\int_0^t\int_{\bf R}\left(v^\varepsilon- v_\infty\right)u_{xxx}^\varepsilon dxds\\[3mm]
 =&\displaystyle\sum\limits_{i=1}^5I_1^i.
 \end{array}
\end{eqnarray}
It is easy to get
\begin{eqnarray}\label{Le-v-firstp2-I11-c}
\arraycolsep=1.5pt
\begin{array}[b]{rl}\displaystyle
 I_1^1\leq&\displaystyle\int_{\bf R} \left(v_0^\varepsilon- v_\infty\right)^2 dx
 +\int_{\bf R}\left(u_{0x}^\varepsilon\right)^2 dx\leq C.
 \end{array}
\end{eqnarray}
Notice that there exists a positive constant $C>0$ such that
\begin{eqnarray}\label{Le-v-firstp3-c}
\left|v^\varepsilon-v_\infty\right|\leq C\left\{
v^\varepsilon\ln\frac{v^\varepsilon}{v_\infty}-\left(v^\varepsilon-v_\infty\right)\right\}\leq
C \eta\left(u^\varepsilon,v^\varepsilon\right),\ \ {\rm for}\ \
v^\varepsilon\in \Omega_1
\end{eqnarray}
and
\begin{eqnarray}\label{Le-v-firstp32-c}
\left|v^\varepsilon-v_\infty\right|^2\leq C\left\{
v^\varepsilon\ln\frac{v^\varepsilon}{v_\infty}-\left(v^\varepsilon-v_\infty\right)\right\}\leq
C \eta\left(u^\varepsilon,v^\varepsilon\right),\ \ {\rm for}\ \
v^\varepsilon\in \Omega_2,
\end{eqnarray}
since
\begin{eqnarray*}
\sup\limits_{v^\varepsilon\in
\Omega_1}\frac{\left|v^\varepsilon-v_\infty\right|}{v^\varepsilon\ln\frac{v^\varepsilon}{v_\infty}
-\left(v^\varepsilon-v_\infty\right)} =\frac{1}{3\ln\frac{3}{2}-1}
\end{eqnarray*}
and
\begin{eqnarray*}
\sup\limits_{v^\varepsilon\in
\Omega_2}\frac{\left|v^\varepsilon-v_\infty\right|^2}{v^\varepsilon\ln\frac{v^\varepsilon}{v_\infty}
-\left(v^\varepsilon-v_\infty\right)}=\frac{v_\infty}{6\ln\frac{3}{2}-2}.
\end{eqnarray*}
Here
\begin{eqnarray*}
\Omega_1=\left\{v^\varepsilon:
\left|v^\varepsilon-v_\infty\right|\geq\frac{v_\infty}{2}\right\}\bigcap
\left\{v^\varepsilon: v^\varepsilon> 0\right\}=\left\{v^\varepsilon:
0< v^\varepsilon\leq\frac{v_\infty}{2}\right\}\cup\left\{
v^\varepsilon\geq\frac{3v_\infty}{2}\right\}
\end{eqnarray*}
 and
\begin{eqnarray*}
\Omega_2=\left\{v^\varepsilon:
\left|v^\varepsilon-v_\infty\right|\leq\frac{v_\infty}{2}\right\}=\left\{v^\varepsilon:
\frac{v_\infty}{2}\leq v^\varepsilon\leq\frac{3v_\infty}{2}\right\}.
\end{eqnarray*}

\noindent Then we have from (2.9) and (2.10)
\begin{eqnarray}\label{Le-v-firstp2-I12-c}
\arraycolsep=1.5pt
\begin{array}[b]{rl}\displaystyle
 I_1^2\leq&\displaystyle
 2\int_{\bf R} \left(v^\varepsilon- v_\infty\right)^2dx
  + \frac{1}{2}\int_{\bf R} \left(u_x^\varepsilon\right)^2 dx\\[3mm]
  \displaystyle
=&\displaystyle
 2\int_{\{x: v^\varepsilon\in \Omega_1\}} \left(v^\varepsilon- v_\infty\right)^2dx
 +2\int_{\{x: v^\varepsilon\in \Omega_2\}} \left(v^\varepsilon- v_\infty\right)^2dx
  + \frac{1}{2}\int_{\bf R} \left(u_x^\varepsilon\right)^2 dx\\[3mm]
\leq&\displaystyle
C\left(1+\left\|v^\varepsilon-v_\infty\right\|_{L^\infty({\bf
R}\times[0,T])}\right)\int_{\bf R} \eta_0 dx
  + \frac{1}{2}\int_{\bf R} \left(u_x^\varepsilon\right)^2 dx.
  \end{array}
\end{eqnarray}
By using integration by part, Cauchy-Schwartz inequality and Lemma
2.1, we have
\begin{eqnarray}\label{Le-v-firstp2-I13-c}
I_1^3=2\int_0^t\int_{\bf R}\left(v_x^\varepsilon\right)^2dxds
\leq\left\|v^\varepsilon\right\|_{L^\infty({\bf R}\times
[0,T])}\int_0^t \int_{\bf R} \frac
{\left(v_x^\varepsilon\right)^2}{v^\varepsilon}dxds,
\end{eqnarray}
\begin{eqnarray}\label{Le-v-firstp2-I14-c}
\arraycolsep=1.5pt
\begin{array}[b]{rl}\displaystyle
 I_1^4=&\displaystyle-4\varepsilon\int_0^t\int_{\bf R}u^\varepsilon u_x^\varepsilon v_x^\varepsilon dxds\\[3mm]
 \leq&\displaystyle \int_0^t\int_{\bf R}v^\varepsilon \left(u_x^\varepsilon\right)^2dxds
 +\left\|u^\varepsilon\right\|_{L^\infty({\bf
R}\times [0,T])}^2\int_0^t \int_{\bf R} \frac
{\left(v_x^\varepsilon\right)^2}{v^\varepsilon}dxds\\[3mm]
\leq&\displaystyle \int_0^t\int_{\bf R}v^\varepsilon
\left(u_x^\varepsilon\right)^2dxds
 +C\left\|u^\varepsilon\right\|_{L^\infty({\bf
R}\times [0,T])}^2
 \end{array}
\end{eqnarray}
and
\begin{eqnarray}\label{Le-v-firstp2-I15-c}
\arraycolsep=1.5pt
\begin{array}[b]{rl}\displaystyle
 I_1^5=&\displaystyle
 2\varepsilon\int_0^t\int_{\bf R} v_{x}^\varepsilon u_{xx}^\varepsilon dxds\\[3mm]
\leq&\displaystyle \varepsilon\int_0^t\int_{\bf R}
\left(u_{xx}^\varepsilon\right)^2 dxds
+\varepsilon\left\|v^\varepsilon\right\|_{L^\infty({\bf R}\times
[0,T])}\int_0^t \int_{\bf R} \frac
{\left(v_x^\varepsilon\right)^2}{v^\varepsilon}dxds.
  \end{array}
\end{eqnarray}
The estimate on $I_2$ is the same as those of $I_1^4$.

{}Finally, we estimate $I_3$. By using Gagliardo-Nirenberg
inequality and Young inequality, it follows that
\begin{eqnarray}\label{Le-v-firstp2-I31-c}
\arraycolsep=1.5pt
\begin{array}[b]{rl}\displaystyle
\left\|u_x^\varepsilon(t)\right\|_{L^3}^3\leq&\displaystyle
C\left\|u_{xx}^\varepsilon(t)\right\|^\frac{7}{4}\left\|u^\varepsilon(t)\right\|^\frac{5}{4}
\leq
\frac{1}{4}\left\|u_{xx}^\varepsilon(t)\right\|^2+C\left\|u^\varepsilon(t)\right\|^{10}.
 \end{array}
\end{eqnarray}
Thus, we have from Lemma 2.1
\begin{eqnarray}\label{Le-v-firstp2-I32-c}
\arraycolsep=1.5pt
\begin{array}[b]{rl}\displaystyle
I_3\leq&\displaystyle \frac{1}{2}\varepsilon\int_0^t\int_{\bf
R}\left(u_{xx}^\varepsilon\right)^2dxds
+C\varepsilon\int_0^t\left\{\int_{\bf R}\left(u^\varepsilon\right)^2dx\right\}^5ds\\[3mm]
\leq&\displaystyle \frac{1}{2}\varepsilon\int_0^t\int_{\bf
R}\left(u_{xx}^\varepsilon\right)^2dxds+C.
 \end{array}
\end{eqnarray}
Substituting (\ref{Le-v-firstp2-I1-c})-(\ref{Le-v-firstp2-I32-c})
into (\ref{Le-v-firstp2-c}), we get that
\begin{eqnarray}\label{Le-v-firstp4-c}
\arraycolsep=1.5pt
\begin{array}[b]{rl}\displaystyle
 &\displaystyle\int_{\bf R}\left(u_x^\varepsilon\right)^2dx
+\int_0^t\int_{\bf R}v^\varepsilon
\left(u_x^\varepsilon\right)^2dxds
 +\varepsilon\int_0^t\int_{\bf R}\left(u_{xx}^\varepsilon\right)^2dxds\\[3mm]
 \leq&\displaystyle C(T)\left(1 + \left\| v^\varepsilon \right\|_{L^\infty({\bf R}\times
[0,T]) } + \left\| u^\varepsilon \right\|_{L^\infty({\bf R}\times
[0,T]) }^2\right).
 \end{array}
\end{eqnarray}
From the Sobolev inequality and Cauchy inequality , we have
\begin{eqnarray}\label{Le-v-firstp5-c}
\left\|u^\varepsilon\right\|_{L^{\infty}({\bf R}\times[0,T])}^2\leq
2\sup\limits_{[0,T]}\left\{\left\|u^\varepsilon(t)\right\|\left\|u_{x}^\varepsilon(t)\right\|\right\}
\leq\sup\limits_{[0,T]}\left\{C\left\|u^\varepsilon(t)\right\|^2+\frac{1}{2C}\left\|u_{x}^\varepsilon(t)\right\|^2\right\}.
\end{eqnarray}
Substituting (\ref{Le-v-firstp5-c}) into $(\ref{Le-v-firstp4-c})$,
using Lemma 2.1, we deduce immediately
\begin{eqnarray}\label{Le-v-firstp6-c}
\arraycolsep=1.5pt
\begin{array}[b]{rl}\displaystyle
 &\displaystyle\int_{\bf R}\left(u_x^\varepsilon\right)^2dx
+\int_0^t\int_{\bf R}v^\varepsilon
\left(u_x^\varepsilon\right)^2dxds
 +\varepsilon\int_0^t\int_{\bf R}\left(u_{xx}^\varepsilon\right)^2dxds\\[3mm]
 \leq&\displaystyle C\left(1 + \left\| v^\varepsilon
\right\|_{L^\infty({\bf R}\times [0,T]) }\right).
 \end{array}
\end{eqnarray}
From Lemma 2.1, (\ref{Le-v-firstp5-c}) and (\ref{Le-v-firstp6-c}),
 we have
\begin{eqnarray}\label{Le-v-firstp61-c}
\left\|u^\varepsilon\right\|_{L^{\infty}({\bf R}\times[0,T])}^2\leq
C\left(1 + \left\| v^\varepsilon \right\|_{L^\infty({\bf R}\times
[0,T]) } \right).
\end{eqnarray}
The rest of the proof of this lemma is similar to those in
\cite{Guo09}. By differentiating the second equation of
(\ref{Equation}) with respect to $x$ once and multiplying the
resulting identity by $2v_x^\varepsilon$, then integrating the final
equation with respect $x$ and $t$ over ${\bf{R}} \times {[0,t]}$, we
have by integrations by parts, (\ref{Le-v-firstp5-c}),
(\ref{Le-v-firstp6-c}) and Lemma 2.1 that
\begin{eqnarray}\label{Le-v-firstp7-c}
\arraycolsep=1.5pt
\begin{array}[b]{rl}\displaystyle
 \displaystyle
\int_{\bf R}\left(v_x^\varepsilon\right)^2dx + \int_0^t\int_{\bf R}
\left(v_{xx}^\varepsilon\right)^2dxds \leq&\displaystyle C\left(1 +
\left\| v^\varepsilon \right\|_{L^\infty({\bf R}\times [0,T])
}^2\right).
\end{array}
\end{eqnarray}
In order to get the $L^\infty({\bf R}\times [0,T])$-norm estimate on
$\left(u^\varepsilon(x,t),v^\varepsilon(x,t)\right)$, take
\begin{eqnarray}\label{Le-v.firstp8-c}
h(v^\varepsilon) = v^\varepsilon\ln \frac{v^\varepsilon}{{v_\infty}}
- \left(v^\varepsilon - v_\infty\right),\ \ \ \ \varphi
(v^\varepsilon) =\left\{\begin{array}{l}
\displaystyle\int_{v_\infty}^{v^\varepsilon}
{\sqrt {h(z)} } dz,\ \ \ \ {\rm for}\ \ v^\varepsilon\geq v_\infty,\\[3mm]
\displaystyle\int_{v^\varepsilon}^{v_\infty} {\sqrt {h(z)} } dz,\ \
\ \ {\rm for}\ \ 0<v^\varepsilon\leq v_\infty.
\end{array}
\right.
\end{eqnarray}
From (2.9) and (2.10), there exists a positive constant $C>0$ such
that
\begin{eqnarray*}
\left|v^\varepsilon-v_\infty\right|\leq
C\left(h(v^\varepsilon)+1\right),
\end{eqnarray*}
which implies
\begin{eqnarray}\label{Le-v.firstp8-1-c}
\left|v^\varepsilon-v_\infty\right|^\frac{1}{2}\leq
C\left(\sqrt{h(v^\varepsilon)}+1\right).
\end{eqnarray}
{}From (\ref{Le-v.firstp8-c}) and (\ref{Le-v.firstp8-1-c}), we have
\begin{eqnarray}\label{Le-v.firstp8-2-c}
\frac{2}{3}\left(v^\varepsilon-v_\infty\right)^\frac{3}{2}
=\int_{v_\infty}^{v^\varepsilon}\left|z-v_\infty\right|^\frac{1}{2}dz\leq
C\int_{v_\infty}^{v^\varepsilon} {\sqrt {h(z)} }
dz+C\left(v^\varepsilon-v_\infty\right),\ \ \ \ {\rm for}\ \
v^\varepsilon\geq v_\infty
\end{eqnarray}
and
\begin{eqnarray}\label{Le-v.firstp8-3-c}
\frac{2}{3}\left(v_\infty-v^\varepsilon\right)^\frac{3}{2}
=\int_{v^\varepsilon}^{v_\infty}\left|z-v_\infty\right|^\frac{1}{2}dz\leq
C\int_{v^\varepsilon}^{v_\infty} {\sqrt {h(z)} }
dz+C\left(v_\infty-v^\varepsilon\right),\ \ \ \ {\rm for}\ \
v^\varepsilon\leq v_\infty.
\end{eqnarray}
On the other hand, we have
\begin{eqnarray}\label{Le-v-firstp9-c-1}
\begin{array}[b]{rl}
\varphi(v^{\varepsilon})= & \displaystyle
\int^x_{-\infty}\varphi(v^{\varepsilon})_ydy =
\int^x_{-\infty}\sqrt{h(v^{\varepsilon})}v^{\varepsilon}_ydy \\
[4mm] \leq & \displaystyle \left\| {\sqrt {h(v^\varepsilon)} }
\right\|\cdot\left\| {v_x^\varepsilon }
  \right\| \leq  C\left(1 +  \left\| v^\varepsilon \right\|_{L^\infty({\bf
R}\times [0,T]) }\right).
 \end{array}
\end{eqnarray}
Thus for any $v^\varepsilon>0$ and some positive constant $C>0$, we
have
\begin{eqnarray}\label{Le-v-firstp9-c}
\begin{array}{rl}
 &\displaystyle C\left(|v^\varepsilon - v_\infty|^\frac{3}{2}-|v^\varepsilon - v_\infty|\right)
 \leq \varphi (v^\varepsilon)
 \leq  C\left(1 +  \left\| v^\varepsilon \right\|_{L^\infty({\bf
R}\times [0,T]) }\right).
\end{array}
\end{eqnarray}
This means that
\begin{eqnarray}\label{Le-v-firstp10-c}
\left\| v^\varepsilon \right\|_{L^\infty ({\bf R}\times [0,T]) }
\leq C,
\end{eqnarray}
and consequently we have from (\ref{Le-v-firstp61-c})
\begin{eqnarray}\label{Le-v-firstp11-c}
\left\| u^\varepsilon \right\|_{L^\infty ({\bf R}\times [0,T]) }
\leq C.
\end{eqnarray}

Combination of (\ref{Le-v-firstp6-c}), (\ref{Le-v-firstp7-c}),
(\ref{Le-v-firstp10-c}) and (\ref{Le-v-firstp11-c}) yields
(\ref{Le-v-first1-c}). This completes the proof of Lemma 2.2.

\noindent\textbf{Remark 2.1.}\ \ The assumption $(u_0^\varepsilon,
v_0^\varepsilon-v_\infty)\in H^2({\bf R})$ implies
\begin{eqnarray*}
\int_{\bf R}\eta_0dx\leq
C\left(\left\|u_0^\varepsilon\right\|^2+\left\|v_0^\varepsilon-v_\infty\right\|^2\right)
\end{eqnarray*}
which and (\ref{Le-v-first1-c}), (\ref{Le-v-firstp2-I12-c}) show
that
\begin{eqnarray*}
\left\|v^\varepsilon(t)-v_\infty\right\|^2\leq
C\left(\left\|u_0^\varepsilon\right\|^2+\left\|v_0^\varepsilon-v_\infty\right\|^2\right).
\end{eqnarray*}

At last, we deduce the $L^\infty([0,T],L^2({\bf R}))$-norm estimate
on $\left(u_{xx}^\varepsilon(x,t), v_{xx}^\varepsilon(x,t)\right)$.
\begin{Lemma} [Second-order energy estimate]
Assume that the assumptions listed in Theorem 1.1 are satisfied.
Then there exists a positive constant $C$, independent of
$\varepsilon$, such that
\begin{eqnarray}\label{Le-v-second1-c}
\begin{array}[b]{rl}
 & \displaystyle\int_{\bf R} \left\{\left(u_{xx}^\varepsilon\right)^2+\left(v_{xx}^\varepsilon\right)^2\right\}dx
+\int_0^t\int_{\bf R}
\left\{\varepsilon\left(u_{xxx}^\varepsilon\right)^2+\left(v_{xxx}^\varepsilon\right)^2\right\}
dxds\leq C.
\end{array}
\end{eqnarray}
\end{Lemma}
\noindent\textbf{Proof.}\ \ Differentiating (\ref{Equation}) with
respect to $x$, multiplying the resulting identity by
$-2u_{xxx}^\varepsilon$ and $-2v_{xxx}^\varepsilon$ respectively,
and integrating the adding result with respect $x$ and $t$ over
${\bf{R}} \times {[0,t]}$, we have from some integrations by parts
that
\begin{eqnarray}\label{Le-v-secondp1-c}
\begin{array}[b]{rl}
 & \displaystyle\int_{\bf R} \left\{\left(u_{xx}^\varepsilon\right)^2+\left(v_{xx}^\varepsilon\right)^2\right\}dx
+\displaystyle2\int_0^t\int_{\bf R} \left\{\varepsilon\left(u_{xxx}^\varepsilon\right)^2
+\left(v_{xxx}^\varepsilon\right)^2\right\} dxds\\[3mm]
 =&\displaystyle
\int_{\bf R}
\left\{\left(u_{0xx}^\varepsilon\right)^2+\left(v_{0xx}^\varepsilon\right)^2\right\}
dx - 2\int_0^t\int_{\bf R} v_{xxx}^\varepsilon \left(u^\varepsilon
v^\varepsilon\right)_{xx}dxds\\[3mm]
&\displaystyle+2\int_0^t\int_{\bf R} v_{xxx}^\varepsilon
u^\varepsilon_{xx}dxds - 2\int_0^t\int_{\bf R} \left\{\varepsilon
\left(u^\varepsilon\right)^2\right\}_{xxx}u_{xx}^\varepsilon dxds\\[3mm]
=&\displaystyle\sum\limits_{i=4}^7I_i.
\end{array}
\end{eqnarray}
The estimates on $I_5$-$I_6$ can be found in \cite{Guo09} as follow:
\begin{eqnarray}\label{Le-v-secondp1-I56-c}
\begin{array}[b]{rl}
 I_5+I_6\leq&\displaystyle
\int_0^t\int_{\bf R} \left(v_{xxx}^\varepsilon\right)^2dxds
+C\left\{1+\int_0^t\int_{\bf R} \left(u_{xx}^\varepsilon\right)^2
dxds\right\}.
\end{array}
\end{eqnarray}
Next we will be devoted to estimating $I_7$.
\begin{eqnarray}\label{Le-v-secondp1-I7-c}
\begin{array}[b]{rl}
 I_7=&\displaystyle
4\varepsilon\int_0^t\int_{\bf R} \left(u_x^\varepsilon\right)^2
u^\varepsilon_{xxx}dxds +4\varepsilon\int_0^t\int_{\bf R}
u^\varepsilon u_{xx}^\varepsilon u^\varepsilon_{xxx}dxds\\[3mm]
=&\displaystyle\sum\limits_{i=1}^2I_7^i.
\end{array}
\end{eqnarray}
{}From Cauchy inequality and Lemma 2.1, one has
\begin{eqnarray}\label{Le-v-secondp1-I71-1-c}
\begin{array}[b]{rl}
 I_7^1\leq&\displaystyle
\frac{1}{4}\varepsilon\int_0^t\int_{\bf R}
\left(u_{xxx}^\varepsilon\right)^2dxds
+16\varepsilon\int_0^t\int_{\bf R}\left(u_x^\varepsilon\right)^4dxds\\[3mm]
\leq&\displaystyle \frac{1}{2}\varepsilon\int_0^t\int_{\bf R}
\left(u_{xxx}^\varepsilon\right)^2dxds +C\int_0^t\left\{\int_{\bf
R}\left(u^\varepsilon\right)^2dx\right\}^7ds\\[3mm]
\leq&\displaystyle \frac{1}{2}\varepsilon\int_0^t\int_{\bf R}
\left(u_{xxx}^\varepsilon\right)^2dxds +C.
\end{array}
\end{eqnarray}
Here we used the fact obtained by Gagliardo-Nirenberg inequality and
Young inequality
\begin{eqnarray}\label{Le-v-secondp1-I71-2-c}
\arraycolsep=1.5pt
\begin{array}[b]{rl}\displaystyle
\left\|u_x^\varepsilon(t)\right\|_{L^4}^4\leq&\displaystyle
C\left\|u_{xxx}^\varepsilon(t)\right\|^\frac{5}{3}\left\|u^\varepsilon(t)\right\|^\frac{7}{3}
\leq
\frac{1}{64}\left\|u_{xxx}^\varepsilon(t)\right\|^2+C\left\|u^\varepsilon(t)\right\|^{14}.
 \end{array}
\end{eqnarray}
In addition, one has by Cauchy inequality and Lemma 2.2
\begin{eqnarray}\label{Le-v-secondp1-I72-c}
\begin{array}[b]{rl}
 I_7^2\leq&\displaystyle
\frac{1}{4}\varepsilon\int_0^t\int_{\bf R}
\left(u_{xxx}^\varepsilon\right)^2dxds
+16\varepsilon\left\|u^\varepsilon\right\|_{L^\infty({\bf
R}\times[0,T])}^2
\int_0^t\int_{\bf R}\left(u_{xx}^\varepsilon\right)^2dxds\\[3mm]
\leq&\displaystyle \frac{1}{4}\varepsilon\int_0^t\int_{\bf R}
\left(u_{xxx}^\varepsilon\right)^2dxds +C\int_0^t\int_{\bf
R}\left(u_{xx}^\varepsilon\right)^2dxds.
\end{array}
\end{eqnarray}
Substituting (\ref{Le-v-secondp1-I56-c})-(\ref{Le-v-secondp1-I72-c})
into (\ref{Le-v-secondp1-c}), we have
\begin{eqnarray}\label{Le-v.secondp2-c}
\begin{array}[b]{rl}
 & \displaystyle\int_{\bf R} \left\{\left(u_{xx}^\varepsilon\right)^2+\left(v_{xx}^\varepsilon\right)^2\right\}dx
+\int_0^t\int_{\bf R} \left\{\varepsilon\left(u_{xxx}^\varepsilon\right)^2
+\left(v_{xxx}^\varepsilon\right)^2\right\} dxds\\[3mm]
 \leq&\displaystyle
C\left\{1+\int_0^t\int_{\bf R} \left(u_{xx}^\varepsilon\right)^2
dxds\right\}.
\end{array}
\end{eqnarray}
(\ref{Le-v-second1-c}) immediately follows by applying the Gronwall
inequality. This completes the proof of Lemma 2.3.

\noindent\textbf{Remark 2.2.}\ \ By Sobolev inequality, Lemmas 2.2
and 2.3, we have
\begin{eqnarray*}
\arraycolsep=1.5pt
\begin{array}[b]{rl}\displaystyle
\left\|u_x^\varepsilon\right\|_{L^\infty({\bf R}\times
[0,T])}\leq&\displaystyle
C\sup\limits_{[0,T]}\left\{\left\|u_x^\varepsilon(t)\right\|^{\frac{1}{2}}\left\|u_{xx}^\varepsilon(t)\right\|^{\frac{1}{2}}\right\}
\leq C
 \end{array}
\end{eqnarray*}
and
\begin{eqnarray*}
\arraycolsep=1.5pt
\begin{array}[b]{rl}\displaystyle
\left\|v_x^\varepsilon\right\|_{L^\infty({\bf R}\times
[0,T])}\leq&\displaystyle
C\sup\limits_{[0,T]}\left\{\left\|v_x^\varepsilon(t)\right\|^{\frac{1}{2}}\left\|v_{xx}^\varepsilon(t)\right\|^{\frac{1}{2}}\right\}
\leq C,
 \end{array}
\end{eqnarray*}
which will be used later.

Combination of Lemmas 2.1-2.3 and Remark 2.1 yields to the
uniform-in-$\varepsilon$ \textit{a priori} estimate
(\ref{Existence3-c}), which implies the global existence of
solutions by combining the local existence and uniqueness. Similar
to \cite{Guo09}, by standard maximal principle, we can prove
(\ref{Existence2-c}) and the priori assumption (\ref{priori
assumption-c}) holds. This completes the proof of global existence
in Theorem 1.1.

\subsection{Convergence rate of zero viscosity limit}
In this subsection, we turn to our another result, which is
concerned with convergence rates of the vanishing diffusion
viscosity. That is, we will give the proof of Theorem 1.2, and it
suffices to show the following Lemma 2.6.

{}First, the global existence of the $H^2$ solutions to the Cauchy
problem (\ref{Limit equation}), (\ref{Limit initial}) was included
in \cite{Guo09}, which is summarized as follows.

\begin{Lemma} [Global existence on the limit problem, cf. \cite{Guo09}]
Assume that the initial data
$\left(u_0^0(x), v_0^0(x)\right)$ satisfies
\begin{eqnarray}\label{Le-l-assum-c}
\inf\limits_{x\in{\bf R}} v_0^0(x)\geq \alpha>0 \quad
{\textrm{and}}\quad \left(u_0^0-\bar{u}, v_0^0-v_\infty\right)\in
H^2({\bf R}).
\end{eqnarray}
Then there exists a unique global solution
$\left(u^0(x,t),v^0(x,t)\right)$ of the Cauchy problem (\ref{Limit
equation}), (\ref{Limit initial}) satisfying
\begin{eqnarray}\label{Le-l-priori-c}
\left(u^0-\bar{u},v^0-v_\infty\right)\in L^\infty
([0,\infty),H^2({\bf R})),\quad v_x^0\in L^2([0,\infty), H^2({\bf
R})).
\end{eqnarray}
Moreover, for any fixed $T>0$, there exists a positive constant
$C(T)>0$ which depends only on $T$ and $\left\|\left(u_0^0-\bar{u},
v_0^0-v_\infty\right)\right\|_{H^2}$ such that
\begin{eqnarray}\label{Le-l-lower-c}
v^0(x,t)\geq C(T).
\end{eqnarray}
Here the smooth monotone function $\bar{u}(x)$ satisfies
$\bar{u}(x)=u_\pm,\ \ \pm x\geq 1$ for fixed constants $u_-$ and
$u_+$.
\end{Lemma}

Next, we need improve the regularity on the solutions to the zero
viscosity limit problem (\ref{Limit equation}), (\ref{Limit
initial}).

\begin{Lemma}[Regularity improved]
Assume that the assumptions listed in Theorem 1.2 are satisfied.
Then there exists a positive constant $C$ such that
\begin{eqnarray}\label{Le-v-reg-c}
\begin{array}[b]{rl}
 & \displaystyle\int_{\bf R} \left\{\left(u_{xxx}^0\right)^2+\left(v_{xxx}^0\right)^2\right\}dx
+\int_0^t\int_{\bf R}\left(v_{xxxx}^0\right)^2dxds\leq C.
\end{array}
\end{eqnarray}
\end{Lemma}
\noindent\textbf{Proof.}\ \ Differentiating (\ref{Limit equation})
with respect to $x$ three times, multiplying the resulting identity
by $2u_{xxx}^0$ and $2v_{xxx}^0$ respectively, and integrating the
adding result with respect $x$ and $t$ over ${\bf{R}} \times
{[0,t]}$, we have from some integrations by parts that
\begin{eqnarray}\label{Le-v-regp1-c}
\begin{array}[b]{rl}
 & \displaystyle\int_{\bf R} \left\{\left(u_{xxx}^0\right)^2+\left(v_{xxx}^0\right)^2\right\}dx
+2\int_0^t\int_{\bf R}\left(v_{xxxx}^0\right)^2dxds\\[3mm]
 =&\displaystyle
\int_{\bf R}
\left\{\left(u_{0xxx}^0\right)^2+\left(v_{0xxx}^0\right)^2\right\}
dx -2\int_0^t\int_{\bf R}u^0 v^0_{xxx}v_{xxxx}^0dxds
-6\int_0^t\int_{\bf R}u_x^0 v^0_{xx}v_{xxxx}^0dxds\\[3mm]
&\displaystyle -6\int_0^t\int_{\bf R}u_{xx}^0 v^0_xv_{xxxx}^0dxds
-2\int_0^t\int_{\bf R}u_{xxx}^0 v^0v_{xxxx}^0dxds
 +2\int_0^t\int_{\bf R}
u^0_{xxx}v_{xxxx}^0
dxds \\[3mm]
=&\displaystyle\sum\limits_{i=1}^6J_i.
\end{array}
\end{eqnarray}
$J_2$-$J_6$ are estimated as follows:
\begin{eqnarray}\label{Le-v-regp1-J2-c}
\begin{array}[b]{rl}
 J_2\leq&\displaystyle
\frac{1}{5}\int_0^t\int_{\bf R}\left(v_{xxxx}^0\right)^2dxds
+C\left\|u^0\right\|_{L^\infty({\bf
R}\times[0,T])}^2\int_0^t\int_{\bf R}\left(v_{xxx}^0\right)^2dxds,
\end{array}
\end{eqnarray}
\begin{eqnarray}\label{Le-v-regp1-J3-c}
\begin{array}[b]{rl}
 J_3\leq&\displaystyle
\frac{1}{5}\int_0^t\int_{\bf R}\left(v_{xxxx}^0\right)^2dxds
+C\left\|u_x^0\right\|_{L^\infty({\bf
R}\times[0,T])}^2\int_0^t\int_{\bf R}\left(v_{xx}^0\right)^2dxds,
\end{array}
\end{eqnarray}
\begin{eqnarray}\label{Le-v-regp1-J4-c}
\begin{array}[b]{rl}
 J_4\leq&\displaystyle
\frac{1}{5}\int_0^t\int_{\bf R}\left(v_{xxxx}^0\right)^2dxds
+C\left\|v_x^0\right\|_{L^\infty({\bf
R}\times[0,T])}^2\int_0^t\int_{\bf R}\left(u_{xx}^0\right)^2dxds,
\end{array}
\end{eqnarray}
\begin{eqnarray}\label{Le-v-regp1-J56-c}
\begin{array}[b]{rl}
 J_5+J_6\leq&\displaystyle
\frac{1}{5}\int_0^t\int_{\bf R}\left(v_{xxxx}^0\right)^2dxds
+C\left(1+\left\|v^0\right\|_{L^\infty({\bf
R}\times[0,T])}^2\right)\int_0^t\int_{\bf
R}\left(u_{xxx}^0\right)^2dxds.
\end{array}
\end{eqnarray}
Substituting (\ref{Le-v-regp1-J2-c})-(\ref{Le-v-regp1-J56-c}) into
(\ref{Le-v-regp1-c}), we get from Lemma 2.4
\begin{eqnarray}\label{Le-v-regp2-c}
\begin{array}[b]{rl}
 & \displaystyle\int_{\bf R} \left\{\left(u_{xxx}^0\right)^2+\left(v_{xxx}^0\right)^2\right\}dx
+\int_0^t\int_{\bf R}\left(v_{xxxx}^0\right)^2dxds\\[3mm]
 \leq&\displaystyle
C\left\{1+\int_0^t\int_{\bf R}\left(u_{xxx}^0\right)^2dxds\right\}.
\end{array}
\end{eqnarray}
(\ref{Le-v-reg-c}) immediately follows by applying the Gronwall
inequality. This completes the proof of Lemma 2.5.

Based on Theorem 1.1 and Lemmas 2.4-2.5, the following
$L^2$-convergence rates can be proved.

\begin{Lemma}[$L^2$-Convergence rates]
Assume that the assumptions listed in Theorem 1.2 are satisfied.
Then there exists a positive constant $C$, independent of
$\varepsilon$, such that
\begin{eqnarray}\label{Le-con1-c}
\arraycolsep=1.5pt
\begin{array}[b]{rl}
\displaystyle\int_{\bf R}\Big[\left(u^\varepsilon-u^0\right)^2
&+\left(v^\varepsilon-v^0\right)^2\Big]dx\\[3mm]
&\displaystyle+\int_0^t\int_{\bf
R}\left[\varepsilon\left(u^\varepsilon-u^0\right)_x^2
+\left(v^\varepsilon-v^0\right)_x^2\right]dxd\tau\leq C\varepsilon^2
\end{array}
\end{eqnarray}
and
\begin{eqnarray}\label{Le-con2-c}
\arraycolsep=1.5pt
\begin{array}[b]{rl}
\displaystyle\int_{\bf R}\Big[\left(u^\varepsilon-u^0\right)_x^2
&+\left(v^\varepsilon-v^0\right)_x^2\Big]dx\\[3mm]
&\displaystyle+\int_0^t\int_{\bf
R}\left[\varepsilon\left(u^\varepsilon-u^0\right)_{xx}^2
+\left(v^\varepsilon-v^0\right)_{xx}^2\right]dxd\tau\leq
C\varepsilon^2.
\end{array}
\end{eqnarray}
\end{Lemma}

\vspace{3mm}\noindent\textbf{Proof.}\ \ Set
\begin{equation}\label{Le-conp1-c}
\psi^\varepsilon=u^\varepsilon-u^0,\ \ \
\theta^\varepsilon=v^\varepsilon-v^0.
\end{equation}
Then we deduce from (\ref{Equation})-(\ref{Initial}) and (\ref{Limit
equation})-(\ref{Limit initial}) that
$\left(\psi^\varepsilon,\theta^\varepsilon\right)(x,t)$ satisfy the
following Cauchy problem:
\begin{equation}\label{Le-conp-equation-c}
\left\{
\begin{array}{l}
\psi_t^\varepsilon+\left(\varepsilon \left(u^\varepsilon\right)^2-\theta^\varepsilon\right)_x
=\varepsilon \psi_{xx}^\varepsilon+\varepsilon u_{xx}^0,\\[3mm]
\theta_t^\varepsilon-\left(\psi^\varepsilon
v^\varepsilon+u^0\theta^\varepsilon\right)_x=\theta_{xx}^\varepsilon,
\end{array}
\right.
\end{equation}
with initial data
\begin{equation}\label{Le-conp-initial-c}
\left(\psi^\varepsilon,
\theta^\varepsilon\right)(x,0)=\left(\psi_0^\varepsilon,
\theta_0^\varepsilon\right).
\end{equation}

\noindent\textbf{\underline{Proof of (\ref{Le-con1-c})}.}

Multiplying the fist and second equation of
(\ref{Le-conp-equation-c}) by $2\psi^\varepsilon$ and
$2\theta^\varepsilon$ respectively, integrating the adding result
with respect $x$ and $t$ over ${\bf R}\times {[0,t]}$, we have
\begin{eqnarray}\label{Le-conp4-c}
\arraycolsep=1.5pt
\begin{array}[b]{rl}\displaystyle
 &\displaystyle\int_{\bf R}\left\{\left(\psi^\varepsilon\right)^2+\left(\theta^\varepsilon\right)^2\right\}dx
+2\int_0^t\int_{\bf R}\left\{\varepsilon\left(\psi_x^\varepsilon\right)^2
+\left(\theta_x^\varepsilon\right)^2\right\}dxds\\[3mm]
 =&\displaystyle \int_{\bf R}\left\{\left(\psi_0^\varepsilon\right)^2+\left(\theta_0^\varepsilon\right)^2\right\}dx
 -4\varepsilon\int_0^t\int_{\bf R}\psi^\varepsilon u^\varepsilon u_x^\varepsilon dxds
+2\varepsilon\int_0^t\int_{\bf R}\psi^\varepsilon
u_{xx}^0dxds\\[3mm]
 &\displaystyle+2\int_0^t\int_{\bf R}\psi^\varepsilon
\theta_x^\varepsilon dxds-2\int_0^t\int_{\bf R}\psi^\varepsilon
v^\varepsilon \theta_x^\varepsilon dxds
 -2\int_0^t\int_{\bf R}u^0 \theta^\varepsilon \theta_x^\varepsilon dxds \\[3mm]
 =&\displaystyle\sum\limits_{i=7}^{12}J_i.
\end{array}
\end{eqnarray}
Here $J_7$-$J_{12}$ are estimated as follows.

{}First, one has
\begin{equation}\label{Le-conp-J7-c}
J_7\leq C\varepsilon^2
\end{equation}
since $\left\|u_0^\varepsilon-u_0^0\right\|_{H^2({\bf
R})}+\left\|v_0^\varepsilon-v_0^0\right\|_{H^2({\bf R})}
=O(\varepsilon)$.

Next, from Cauchy-Schwarz inequality, Theorem 1.1 and Lemma 2.4 with
$\bar{u}(x)\equiv 0$, we obtain
\begin{eqnarray}\label{Le-conp-J8-c}
\arraycolsep=1.5pt
\begin{array}[b]{rl}\displaystyle
J_8 \leq&\displaystyle
4\varepsilon^2\left\|u^\varepsilon\right\|_{L^\infty({\bf R }\times
[0,T])}^2\int_0^t\int_{\bf R}\left(u_x^\varepsilon\right)^2dxds
+\int_0^t\int_{\bf R} \left(\psi^\varepsilon\right)^2dxds\\[3mm]
\leq&\displaystyle C\varepsilon^2+\int_0^t\int_{\bf R}
\left(\psi^\varepsilon\right)^2dxds,
\end{array}
\end{eqnarray}
\begin{eqnarray}\label{Le-conp-J9-c}
\arraycolsep=1.5pt
\begin{array}[b]{rl}\displaystyle
J_9 \leq&\displaystyle 4\varepsilon^2\int_0^t\int_{\bf
R}\left(u_{xx}^0\right)^2dxds
+\int_0^t\int_{\bf R} \left(\psi^\varepsilon\right)^2dxds\\[3mm]
\leq&\displaystyle C\varepsilon^2+\int_0^t\int_{\bf R}
\left(\psi^\varepsilon\right)^2dxds
\end{array}
\end{eqnarray}
and
\begin{eqnarray}\label{Le-conp-J10-J12-c}
\arraycolsep=1.5pt
\begin{array}[b]{rl}\displaystyle
\sum\limits_{i=10}^{12}J_i \leq&\displaystyle
2\left(1+\left\|v^\varepsilon\right\|_{L^\infty({\bf R }\times
[0,T])}^2+\left\|u^0\right\|_{L^\infty({\bf R }\times
[0,T])}^2\right)\int_0^t\int_{\bf R} \left\{\left(\psi^\varepsilon\right)^2
+\left(\theta^\varepsilon\right)^2\right\}dxds\\[3mm]
&\displaystyle+\frac{3}{2}\int_0^t\int_{\bf R}\left(\theta_x^\varepsilon\right)^2dxds \\[3mm]
\leq&\displaystyle C\int_0^t\int_{\bf R}
\left\{\left(\psi^\varepsilon\right)^2+\left(\theta^\varepsilon\right)^2\right\}dxds
+\frac{3}{2}\int_0^t\int_{\bf
R}\left(\theta_x^\varepsilon\right)^2dxds.
\end{array}
\end{eqnarray}
Substituting (\ref{Le-conp-J7-c})-(\ref{Le-conp-J10-J12-c}) into
(\ref{Le-conp4-c}), we get
\begin{eqnarray}\label{Le-conp8-c}
\arraycolsep=1.5pt
\begin{array}[b]{rl}\displaystyle
 &\displaystyle\int_{\bf R}\left\{\left(\psi^\varepsilon\right)^2+\left(\theta^\varepsilon\right)^2\right\}dx
+\int_0^t\int_{\bf R}\left\{\varepsilon\left(\psi_x^\varepsilon\right)^2
+\left(\theta_x^\varepsilon\right)^2\right\}dxds\\[3mm]
 \leq&\displaystyle C \varepsilon^2+C\int_0^t\int_{\bf R}
\left\{\left(\psi^\varepsilon\right)^2+\left(\theta^\varepsilon\right)^2\right\}dxds.
\end{array}
\end{eqnarray}
We immediately get (\ref{Le-con1-c}) from Gronwall inequality.

\noindent\textbf{\underline{Proof of (\ref{Le-con2-c})}.}

Multiplying the fist and second equation of
(\ref{Le-conp-equation-c}) by $-2\psi_{xx}^\varepsilon$ and
$-2\theta_{xx}^\varepsilon$ respectively, integrating the adding
result with respect $x$ and $t$ over ${\bf R}\times {[0,t]}$, we
have
\begin{eqnarray}\label{Le-conp9-c}
\arraycolsep=1.5pt
\begin{array}[b]{rl}\displaystyle
 &\displaystyle\int_{\bf R}\left\{\left(\psi_x^\varepsilon\right)^2+\left(\theta_x^\varepsilon\right)^2\right\}dx
+2\int_0^t\int_{\bf R}\left\{\varepsilon\left(\psi_{xx}^\varepsilon\right)^2
+\left(\theta_{xx}^\varepsilon\right)^2\right\}dxds\\[3mm]
 =&\displaystyle \int_{\bf R}\left\{\left(\psi_{0x}^\varepsilon\right)^2
 +\left(\theta_{0x}^\varepsilon\right)^2\right\}dx
 +4\varepsilon\int_0^t\int_{\bf R}\psi_{xx}^\varepsilon u^\varepsilon u_x^\varepsilon dxds
-2\varepsilon\int_0^t\int_{\bf R}\psi_{xx}^\varepsilon
u_{xx}^0dxds\\[3mm]
 &\displaystyle-2\int_0^t\int_{\bf R}\psi_{xx}^\varepsilon
\theta_x^\varepsilon dxds
-2\int_0^t\int_{\bf R}\left(\psi^\varepsilon v^\varepsilon
+u^0 \theta^\varepsilon\right)_x \theta_{xx}^\varepsilon dxds\\[3mm]
 =&\displaystyle\sum\limits_{i=13}^{17}J_i.
\end{array}
\end{eqnarray}
Here $J_{13}$-$J_{17}$ are estimated as follows.

{}First, similar to the estimate (\ref{Le-conp-J7-c}), one has
\begin{equation}\label{Le-conp-J13-c}
J_{13}\leq C\varepsilon^2.
\end{equation}

Next, from Cauchy-Schwarz inequality, Theorem 1.1, Lemmas 2.4-2.5,
Remark 2.2 and (\ref{Le-con1-c}), we obtain
\begin{eqnarray}\label{Le-conp-J14-c}
\arraycolsep=1.5pt
\begin{array}[b]{rl}\displaystyle
J_{14}=&\displaystyle -4\varepsilon\int_0^t\int_{\bf
R}\psi_x^\varepsilon \left(u_x^\varepsilon\right)^2 dxds
-4\varepsilon\int_0^t\int_{\bf R}\psi_x^\varepsilon u^\varepsilon u_{xx}^\varepsilon dxds\\[3mm]
\leq&\displaystyle
C\varepsilon^2\left\|u_x^\varepsilon\right\|_{L^\infty({\bf R
}\times [0,T])}^2\int_0^t\int_{\bf
R}\left(u_x^\varepsilon\right)^2dxds
+C\varepsilon^2\left\|u^\varepsilon\right\|^2_{L^\infty({\bf R
}\times [0,T])}\int_0^t\int_{\bf
R}\left(u_{xx}^\varepsilon\right)^2dxds\\[3mm]
&\displaystyle+\int_0^t\int_{\bf R} \left(\psi_x^\varepsilon\right)^2dxds\\[3mm]
\leq&\displaystyle C\varepsilon^2+\int_0^t\int_{\bf R}
\left(\psi_x^\varepsilon\right)^2dxds,
\end{array}
\end{eqnarray}
\begin{eqnarray}\label{Le-conp-J15-c}
\arraycolsep=1.5pt
\begin{array}[b]{rl}\displaystyle
J_{15} \leq&\displaystyle \varepsilon^2\int_0^t\int_{\bf
R}\left(u_{xxx}^0\right)^2dxds
+\int_0^t\int_{\bf R} \left(\psi_x^\varepsilon\right)^2dxds\\[3mm]
\leq&\displaystyle C\varepsilon^2+\int_0^t\int_{\bf R}
\left(\psi_x^\varepsilon\right)^2dxds,
\end{array}
\end{eqnarray}
\begin{eqnarray}\label{Le-conp-J16-c}
\arraycolsep=1.5pt
\begin{array}[b]{rl}\displaystyle
J_{16}=&\displaystyle 2\int_0^t\int_{\bf R}\psi_x^\varepsilon
\theta_{xx}^\varepsilon dxds\\[3mm]
\leq&\displaystyle \int_0^t\int_{\bf
R}\left(\psi_x^\varepsilon\right)^2dxds
+\frac{1}{4}\int_0^t\int_{\bf R}
\left(\theta_{xx}^\varepsilon\right)^2dxds
\end{array}
\end{eqnarray}
and
\begin{eqnarray}\label{Le-conp-J17-c}
\arraycolsep=1.5pt
\begin{array}[b]{rl}\displaystyle
 J_{17}=&\displaystyle-2\int_0^t\int_{\bf R}\psi_x^\varepsilon v^\varepsilon\theta_{xx}^\varepsilon dxds
 -2\int_0^t\int_{\bf R}\psi^\varepsilon v_x^\varepsilon\theta_{xx}^\varepsilon dxds
 -2\int_0^t\int_{\bf R}u_x^0 \theta^\varepsilon\theta_{xx}^\varepsilon
 dxds\\[3mm]
 &\displaystyle-2\int_0^t\int_{\bf R}u^0 \theta_x^\varepsilon\theta_{xx}^\varepsilon dxds\\[3mm]
 =&\displaystyle\sum\limits_{i=1}^4J_{17}^i.
\end{array}
\end{eqnarray}
Here
\begin{eqnarray}\label{Le-conp-J171-c}
\arraycolsep=1.5pt
\begin{array}[b]{rl}\displaystyle
J_{17}^1 \leq&\displaystyle \frac{1}{4}\int_0^t\int_{\bf
R}\left(\theta_{xx}^\varepsilon\right)^2dxds
+C\left\|v^\varepsilon\right\|_{L^\infty({\bf R }\times
[0,T])}^2\int_0^t\int_{\bf
R}\left(\psi_x^\varepsilon\right)^2dxds\\[3mm]
\leq&\displaystyle \frac{1}{4}\int_0^t\int_{\bf
R}\left(\theta_{xx}^\varepsilon\right)^2dxds+C\int_0^t\int_{\bf R}
\left(\psi_x^\varepsilon\right)^2dxds,
\end{array}
\end{eqnarray}
\begin{eqnarray}\label{Le-conp-J172-c}
\arraycolsep=1.5pt
\begin{array}[b]{rl}\displaystyle
J_{17}^2 \leq&\displaystyle \frac{1}{4}\int_0^t\int_{\bf
R}\left(\theta_{xx}^\varepsilon\right)^2dxds
+C\left\|v_x^\varepsilon\right\|_{L^\infty({\bf R }\times
[0,T])}^2\int_0^t\int_{\bf
R}\left(\psi^\varepsilon\right)^2dxds\\[3mm]
\leq&\displaystyle \frac{1}{4}\int_0^t\int_{\bf
R}\left(\theta_{xx}^\varepsilon\right)^2dxds+C\varepsilon^2,
\end{array}
\end{eqnarray}
\begin{eqnarray}\label{Le-conp-J173-4-c}
\arraycolsep=1.5pt
\begin{array}[b]{rl}\displaystyle
J_{17}^3+J_{17}^4 \leq&\displaystyle \frac{1}{4}\int_0^t\int_{\bf
R}\left(\theta_{xx}^\varepsilon\right)^2dxds
+C\left\|u_x^0\right\|_{L^\infty({\bf R }\times
[0,T])}^2\int_0^t\int_{\bf
R}\left(\theta^\varepsilon\right)^2dxds\\[3mm]
&\displaystyle+C\left\|u^0\right\|_{L^\infty({\bf R }\times
[0,T])}^2\int_0^t\int_{\bf
R}\left(\theta_x^\varepsilon\right)^2dxds\\[3mm]
\leq&\displaystyle \frac{1}{4}\int_0^t\int_{\bf
R}\left(\theta_{xx}^\varepsilon\right)^2dxds+C\varepsilon^2.
\end{array}
\end{eqnarray}
Substituting (\ref{Le-conp-J13-c})-(\ref{Le-conp-J173-4-c}) into
(\ref{Le-conp9-c}), we get
\begin{eqnarray}\label{Le-conp16-c}
\arraycolsep=1.5pt
\begin{array}[b]{rl}\displaystyle
 &\displaystyle\int_{\bf R}\left\{\left(\psi_x^\varepsilon\right)^2+\left(\theta_x^\varepsilon\right)^2\right\}dx
+\int_0^t\int_{\bf R}\left\{\varepsilon\left(\psi_{xx}^\varepsilon\right)^2+\left(\theta_{xx}^\varepsilon\right)^2\right\}dxds\\[3mm]
 \leq&\displaystyle C \varepsilon^2+C\int_0^t\int_{\bf R}
\left\{\left(\psi_x^\varepsilon\right)^2+\left(\theta_x^\varepsilon\right)^2\right\}dxds.
\end{array}
\end{eqnarray}
We immediately get (\ref{Le-con2-c}) from Gronwall inequality. This
completes the proof of Lemma 2.6.

{}Finally, based on Lemma 2.6,  we can prove Theorem 1.2. In fact,
using Sobolev inequality, we also have from (\ref{Le-con1-c}) and
(\ref{Le-con2-c})
\begin{eqnarray}\label{Th-conp1-c}
\arraycolsep=1.5pt
\begin{array}[b]{rl}\displaystyle
\left\|\left(u^\varepsilon-u^0\right)(t)\right\|_{L^\infty({\bf R})}
\leq&\displaystyle
C\left\|\left(u^\varepsilon-u^0\right)(t)\right\|^\frac{1}{2}
\left\|\left(u^\varepsilon-u^0\right)_x(t)\right\|^\frac{1}{2}\\[3mm]
\leq&\displaystyle C\varepsilon
\end{array}
\end{eqnarray}
and
\begin{eqnarray}\label{Th-conp2-c}
\arraycolsep=1.5pt
\begin{array}[b]{rl}\displaystyle
\left\|\left(v^\varepsilon-v^0\right)(t)\right\|_{L^\infty({\bf R})}
\leq&\displaystyle
C\left\|\left(v^\varepsilon-v^0\right)(t)\right\|^\frac{1}{2}
\left\|\left(v^\varepsilon-v^0\right)_x(t)\right\|^\frac{1}{2}\\[3mm]
\leq&\displaystyle C\varepsilon.
\end{array}
\end{eqnarray}
This completes the proof of Theorem 1.2.
\section{Initial-boundary value problem}
\setcounter{equation}{0} In this section, we are concerned with the
global existence of large-amplitude $H^2$ solutions to the
initial-boundary value problem (\ref{Equation}), (\ref{Initial1})
and (\ref{boundary}). In particular, we investigate the convergence
rates as the diffusion parameter $\varepsilon$ goes to zero.

{}By the method similar to those in Subsection 2.1, we can prove
(\ref{Existence1-cb}) and (\ref{Existence3-cb}) in Theorem 1.3 ({\rm
i}). Now we devoted ourselves to claiming (\ref{Existence2-cb}) in
Theorem 1.3 ({\rm i}).

\noindent\textbf{\underline{Proof of (\ref{Existence2-cb}).}}\ \
From Sobolev inequality and (\ref{Existence3-cb}), one has
$$
\begin{array} {rl}
&\displaystyle\left\| {u_x^\varepsilon }
\right\|_{L^\infty([0,1]\times [0,T]) } \leq C.
\end{array}
$$
Now set
$$
M=\mathop {\sup }\limits_{[0,1] \times [0,T]}\Big\{ \left|
{u_x^\varepsilon (x,t)} \right|\Big\},\ \ \
w^\varepsilon(x,t)=e^{-\lambda t}v^\varepsilon(x,t),
$$
where $\lambda>M$ is a positive constant.

We can deduce from (\ref{Equation}), (\ref{Initial1}) and
(\ref{boundary}) that $w^\varepsilon(x,t)$ satisfies the following
initial-boundary value problem
\begin{equation}\label{IBVP-w}
\left\{
\begin{array}{l}
w_t^\varepsilon  - u^\varepsilon w_x^\varepsilon
+\left(\lambda- u_x^\varepsilon \right)w^\varepsilon = w_{xx}^\varepsilon,\ \ x\in(0,1),\quad 0<t<T,\\[3mm]
w^\varepsilon(x,0)=v_0^\varepsilon(x)\geq \alpha>0,\ \ x\in[0,1],\\[3mm]
w_x^\varepsilon(0,t)=w_x^\varepsilon(1,t)=0,\ \ 0\leq t\leq T.
\end{array}
\right.
\end{equation}
The standard maximal principle tells us that
\begin{eqnarray}\label{Le-l.lowerp1}
w^\varepsilon(x,t)\geq 0,\quad (x,t)\in [0,1]\times [0,T].
\end{eqnarray}
In fact, if $w^\varepsilon(x,t)$ attains its negative minimum at
some point $(x_0,t_0)\in [0,1]\times(0,T]$. Thus
$w^\varepsilon(x,t)$ satisfies
\begin{eqnarray*}
w_t^\varepsilon(x_0,t_0)\leq 0,\ \ w_x^\varepsilon(x_0,t_0)=0,\ \
w^\varepsilon(x_0,t_0)<0,\ \ w_{xx}^\varepsilon(x_0,t_0)\geq 0,
\end{eqnarray*}
which contradict to the first equation of (\ref{IBVP-w}) since
$\lambda-u_x^\varepsilon>0$.

One can deduce that $w^\varepsilon(x,t)$ satisfies from
(\ref{IBVP-w}) and (\ref{Le-l.lowerp1})
\begin{eqnarray}\label{Le-l.lowerp2}
w_t^\varepsilon(x,t)+(\lambda+M)w^\varepsilon(x,t)-u^\varepsilon(x,t)w_{x}^\varepsilon(x,t)\geq
w_{xx}^\varepsilon(x,t).
\end{eqnarray}
Let $H^\varepsilon(x,t)=e^{-\alpha t}\left(e^{(\lambda+M)t}w^\varepsilon(x,t)-\alpha\right)$, we have
from (\ref{Le-l.lowerp2}) that
$$
\left\{\begin{array}{rl}
&H_t^\varepsilon(x,t)+(\lambda+M) H^\varepsilon(x,t) -u^\varepsilon(x,t)H_x^\varepsilon(x,t)
\geq H_{xx}^\varepsilon(x,t),\ \ x\in (0,1),\quad 0< t< T,\\[3mm]
&H^\varepsilon(x,0)=w^\varepsilon(x,0)-\alpha\geq 0,\ \ x\in [0,1],\\[3mm]
&H_x^\varepsilon(0,t)=H_x^\varepsilon(1,t)=0,\ \ 0\leq t\leq T.
\end{array}
\right.
$$
By the standard maximal principle, we have
$$
H^\varepsilon(x,t)\geq 0, \quad (x,t)\in [0,1]\times [0,T].
$$
Thus for any $(x,t)\in [0,1]\times {[0,T]}$, $v^\varepsilon(x,t)$
satisfies
\begin{eqnarray}\label{3.28}
v^\varepsilon(x,t)\geq \alpha e^{-Mt}\geq \alpha e^{-MT}.
\end{eqnarray}
This completes the proof of (\ref{Existence2-cb}).

Similar results to those of Theorem 1.3 ({\rm i}) can be obtained on
the initial-boundary value problem on (\ref{Limit
equation1})-(\ref{boundary11}) as follows.
\begin{Lemma}
Let $\left(u_0^0,v_0^0\right)\in H^2({[0,1]})$ and assume that there
exists a positive constant $\alpha>0$ such that
$\inf\limits_{x\in{[0,1]}}v_0^0(x)\geq \alpha$. Then there exists a
unique global solution $\left(u^0(x,t),v^0(x,t)\right)$ of the
initial-boundary value problem (\ref{Limit
equation1})-(\ref{boundary11}) satisfying
\begin{equation}\label{Le-l-priori-c1}
\left\{
\begin{array}{l}
\left(u^0, v^0\right) \in
L^\infty([0,\infty),H^2([0,1])),\\[3mm]
 \frac{v_x^0}{\sqrt{v^0}}\in L^2([0,\infty),L^2([0,1])),\quad
 v_{xx}^0\in L^2([0,\infty),H^1([0,1])).
 \end{array}
 \right.
\end{equation}
Moreover, for any fixed $T>0$, there exists a positive constant
$C(T)>0$ which depends only on $T$ and $\left\|\left(u_0^0,
v_0^0\right)\right\|_{H^2}$ such that
\begin{eqnarray}\label{Le-l-lower-c1}
v^0(x,t)\geq C(T).
\end{eqnarray}
\end{Lemma}
Based on Theorem 1.3 ({\rm i}) and Lemma 3.1, the following
$L^2$-convergence rates can be proved.
\begin{Lemma}[$L^2$-Convergence rates]
Assume that the assumptions listed in Theorem 1.3 are satisfied.
Then there exists a positive constant $C$, independent of
$\varepsilon$, such that
\begin{eqnarray}\label{Le-con1-cb}
\arraycolsep=1.5pt
\begin{array}[b]{rl}
\displaystyle\int_0^1\Big[\left(u^\varepsilon-u^0\right)^2
&+\left(v^\varepsilon-v^0\right)^2\Big]dx\\[3mm]
&\displaystyle+\int_0^t\int_0^1\left[\varepsilon\left(u^\varepsilon-u^0\right)_x^2
+\left(v^\varepsilon-v^0\right)_x^2\right]dxd\tau\leq C\varepsilon^2
\end{array}
\end{eqnarray}
and
\begin{eqnarray}\label{Le-con2-cb}
\arraycolsep=1.5pt
\begin{array}[b]{rl}
\displaystyle\int_0^1\Big[\left(u^\varepsilon-u^0\right)_x^2
&+\left(v^\varepsilon-v^0\right)_x^2\Big]dx\\[3mm]
&\displaystyle+\int_0^t\int_0^1\left[\varepsilon\left(u^\varepsilon-u^0\right)_{xx}^2
+\left(v^\varepsilon-v^0\right)_{xx}^2\right]dxd\tau\leq
C\varepsilon.
\end{array}
\end{eqnarray}
\end{Lemma}

\vspace{3mm}\noindent\textbf{Proof.}\ \ Set
\begin{equation}\label{Le-conp1-c}
\psi^\varepsilon=u^\varepsilon-u^0,\ \ \
\theta^\varepsilon=v^\varepsilon-v^0.
\end{equation}
Then we deduce from (\ref{Equation}), (\ref{Initial1}),
(\ref{boundary}) and (\ref{Limit equation1})-(\ref{boundary11}) that
$\left(\psi^\varepsilon,\theta^\varepsilon\right)(x,t)$ satisfy the
following initial-boundary value problem:
\begin{equation}\label{Le-conp-equation-c3}
\left\{
\begin{array}{l}
\psi_t^\varepsilon+\left(\varepsilon
\left(u^\varepsilon\right)^2-\theta^\varepsilon\right)_x
=\varepsilon \psi_{xx}^\varepsilon+\varepsilon u_{xx}^0,\\[3mm]
\theta_t^\varepsilon-\left(\psi^\varepsilon
v^\varepsilon+u^0\theta^\varepsilon\right)_x=\theta_{xx}^\varepsilon,
\end{array}
\right.
\end{equation}
with initial data
\begin{equation}\label{Le-conp-initial-c}
\left(\psi^\varepsilon,
\theta^\varepsilon\right)(x,0)=\left(\psi_0^\varepsilon,
\theta_0^\varepsilon\right),
\end{equation}
and the boundary conditions
\begin{equation}\label{boundary-c}
\left(\psi^\varepsilon,
\theta^\varepsilon_x\right)(0,t)=\left(\psi^\varepsilon,\theta_x^\varepsilon\right)(1,t)=(0,0),\
\ \ t\geq 0.
\end{equation}
By the similar method to those in Subsection 2.2, we can obtain
(\ref{Le-con1-cb}) and (\ref{Le-con2-cb}). We note that the
convergence rate in (\ref{Le-con1-cb}) is the same as that shown in
(\ref{Le-con1-c}). However, the rate in (\ref{Le-con2-cb}) is slower
than one obtained in (\ref{Le-con2-c}), which is caused by the
boundary effect. In fact, unlike Lemma 2.5, the regularity on
solutions to the initial-boundary value problem (\ref{Limit
equation1})-(\ref{boundary11}) can not be improved. Consequently,
the term $J_{15}=-2\varepsilon\int_0^t\int_0^1\psi_{xx}^\varepsilon
u_{xx}^0dxds $ (see (2.59) above) can only be controlled as follows:
\begin{eqnarray}\label{Le-conp-J15-c3}
\arraycolsep=1.5pt
\begin{array}[b]{rl}\displaystyle
-2\varepsilon\int_0^t\int_0^1\psi_{xx}^\varepsilon u_{xx}^0dxds
\leq&\displaystyle
C\varepsilon\int_0^t\int_0^1\left(u_{xx}^0\right)^2dxds
+\varepsilon\int_0^t\int_0^1 \left(\psi_{xx}^\varepsilon\right)^2dxds\\[3mm]
\leq&\displaystyle C\varepsilon+\varepsilon\int_0^t\int_0^1
\left(\psi_{xx}^\varepsilon\right)^2dxds.
\end{array}
\end{eqnarray}
{}Finally, based on Lemma 3.2,  we can prove Theorem 1.3 (ii). In
fact, using Sobolev inequality, we also have from (\ref{Le-con1-cb})
and (\ref{Le-con2-cb})
\begin{eqnarray}\label{Th-conp1-c3}
\arraycolsep=1.5pt
\begin{array}[b]{rl}\displaystyle
\left\|\left(u^\varepsilon-u^0\right)(t)\right\|_{L^\infty[0,1]}
\leq&\displaystyle
C\left\|\left(u^\varepsilon-u^0\right)(t)\right\|^\frac{1}{2}
\left\|\left(u^\varepsilon-u^0\right)_x(t)\right\|^\frac{1}{2}\\[3mm]
\leq&\displaystyle C\varepsilon^{\frac{3}{4}}
\end{array}
\end{eqnarray}
and
\begin{eqnarray}\label{Th-conp2-c3}
\arraycolsep=1.5pt
\begin{array}[b]{rl}\displaystyle
\left\|\left(v^\varepsilon-v^0\right)(t)\right\|_{L^\infty[0,1]}
\leq&\displaystyle
C\left\|\left(v^\varepsilon-v^0\right)(t)\right\|^\frac{1}{2}
\left\|\left(v^\varepsilon-v^0\right)_x(t)\right\|^\frac{1}{2}\\[3mm]
\leq&\displaystyle C\varepsilon^{\frac{3}{4}}.
\end{array}
\end{eqnarray}
This completes the proof of Theorem 1.3.

\section{Appendix: Derivation of models}
\setcounter{equation}{0} In this section, we give the derives of
conservation laws (\ref{Equation}) and (\ref{Limit equation}) for
the convenience of the readers respectively, cf. \cite{Li11,
Zhang07, Guo09}.

As in \cite{Li11}, one considers the case where $f(c)=\alpha c$
$(\alpha>0)$ in the Keller-Segel model (\ref{Keller-Segel}), which
means that oxygen is consumed only when cells (bacteria) encounter
the chemical (oxygen). The crucial step in the derivation is to make
a change of variable through the Hopf-Cole transformation
\begin{equation}\label{Hopf-Cole}
\begin{array}{l}
v=-c^{-1}c_x=-(\ln c)_x
\end{array}
\end{equation}
which was first introduced in \cite{Wang08} for a chemotaxis model
proposed in \cite{Levine97} describing the chemotactic movement for
non-diffusible chemicals (i.e. $\varepsilon=0$), and was later
applied in \cite{Li09, Li10} to study the nonlinear stability of
traveling wave solutions. It turns out this transformation also
extends its capacity to the full Keller-Segel model
(\ref{Keller-Segel}) for $\varepsilon\neq 0$. Indeed, with the
Hopf-Cole transformation (\ref{Hopf-Cole}), they derive the
following viscous conservation laws from (\ref{Keller-Segel})
\begin{equation}\label{viscous laws}
\left\{
\begin{array}{l}
u_t+\left(\varepsilon u^2-\alpha v\right)_x=\varepsilon u_{xx},\\[3mm]
v_t-\chi\left(u v\right)_x=Dv_{xx}.
\end{array}
\right.
\end{equation}
Substituting the scalings
\begin{equation}\label{scalings-1}
\begin{array}{l}
\tilde{t}=\alpha t,\ \ \tilde{x}=\sqrt{\frac{\alpha}{\chi}}x,\ \
\tilde{v}=\sqrt{\frac{\alpha}{\chi}}v,\ \ \tilde{D}=\frac{D}{\chi},\
\ \tilde{\varepsilon}=\frac{\varepsilon}{\chi}
\end{array}
\end{equation}
into (\ref{viscous laws}) and dropping the tildes for convenience,
we obtain the system of conservation laws (\ref{Equation}). By the
transformation mentioned above, we get conservation laws (\ref{Limit
equation}) from (\ref{Othmer-Stevens}).

\vskip 0.6cm\noindent\textbf{Acknowledgement:}\ \ The second author
was supported by the Natural Science Foundation of China (The Youth
Foundation) $\#$10901068 and CCNU Project (No. CCNU09A01004) and the
third author was supported by the National Natural Science
Foundation of China $\#$10625105, $\#$11071093, the PhD specialized
grant of the Ministry of Education of China $\#$20100144110001, and
the Special Fund for Basic Scientific Research  of Central Colleges
$\#$CCNU10C01001.


\begin{thebibliography}{99}
\bibitem{Adler66} J. Adler, Chemotaxis in bacteria, \textit{Science}, \textbf{153} (1966), 708-716.

\bibitem{Adler69} J. Adler, Chemoreceptors in bacteria, \textit{Science}, \textbf{166} (1969), 1588-1597.

\bibitem{Chen08} K.M. Chen and C.J. Zhu, The zero diffusion limit
for nonlinear hyperbolic system with damping and diffusion,
\textit{J. Hyperbolic Differ. Equ.}, \textbf{5} (2008), 767-783.

\bibitem{Frid99} H. Frid and V. Shelukhin, Boundary layers for the
Navier-Stokes equations of compressible fluids, \textit{Comm. Math.
Phys.}, \textbf{208} (1999), 309-330.

\bibitem{Guo09} J. Guo, J.X. Xiao, H.J. Zhao and C.J. Zhu, Global solutions
to a hyperbolic-parabolic coupled system with large initial data,
 \textit{Acta Math. Sci. Ser. B}, \textbf{29} (2009), 629-641.

\bibitem{Hillen04} T. Hillen and A. Potapov, The one-dimensional chemotaxis
model:\ global existence and asymptotic profile, \textit{Math.
Methods Appl. Sci.}, \textbf{27} (2004), 1783-1801.

\bibitem{Jiang09} S. Jiang and J.W. Zhang, Boundary layers for the
Navier-Stokes equations of compressible heat-conducting flows with
cylindrical symmetry, \textit{SIAM J. Math. Anal.}, \textbf{41}
(2009), 237-268.

\bibitem{Keller71} E.F. Keller and L.A. Segel, Traveling bands of chemotactic bacteria: a
theoretical analysis, \textit{J. Theoret. Biol.}, \textbf{26}
(1971), 235-248.

\bibitem{Levine97} H.A. Levine and B.D. Sleeman, A system of reaction
diffusion equtions arising in the theory of reinforced random walks,
\textit{SIAM J. Appl. Math.}, \textbf{57} (1997), 683-730.

\bibitem{Levine01} H.A. Levine, B.D. Sleeman and M.N. Hamilton, Mathematical
modeling of the onset of capillary formation initating angiogenesis,
\textit{J. Math. Biol.}, \textbf{42} (2001), 195-238.

\bibitem{Li111} D. Li, T. Li and K. Zhao, On a hyperbolic-parabolic
system modeling chemotaxis, \textit{Math. Models Methods Appl.
Sci.}, \textbf{1} (2011), 1-21.

\bibitem{Li09} T. Li and Z.A. Wang, Nonlinear stability of traveling waves
to a hyperbolic-parabolic system modeling chemotaxis, \textit{SIAM
J. Appl. Math.}, \textbf{70} (2009), 1522-1541.

\bibitem{Li10} T. Li and Z.A. Wang, Nonlinear stability of large amplitude viscous shock waves
of a generalized hyperbolic-parabolic system arising in chemotaxis,
\textit{Math. Models Methods Appl. Sci.}, \textbf{20} (2010),
1967-1998.

\bibitem{Li11} T. Li and Z.A. Wang, Asymptotic nonlinear stability of traveling waves to
conservation laws arising from chemotaxis, \textit{J. Differential
Equations}, \textbf{250} (2011), 1310-1333.

\bibitem{Nagai91} T. Nagai and T. Ikeda, Traveling waves in a chemotaxis model, \textit{J. Math.
Biol.}, \textbf{30} (1991), 169-184.


\bibitem{Othmer97} H.G. Othmer and A. Stevens, Aggregation, blowup, and
collapse: the ABC's of taxis in reinforced random walks,
\textit{SIAM J. Appl. Math.}, \textbf{57} (1997), 1044-1081.

\bibitem{Painter11} K.J. Painter and T. Hillen, Spatio-temporal chaos in a chemotaxis model,
\textit{Phys. D}, \textbf{240} (2011), 363-375.

\bibitem{Ruan11} L.Z. Ruan and C.J. Zhu, Boundary layer for nonlinear evolution equations with damping and diffusion,
\textit{Discrete Contin. Dyn. Syst. Ser. A}, (2011).

\bibitem{Sleeman01} B.D. Sleeman and H.A. Levine, Partial differential
equations of chemotaxis and angiogenesis, \textit{Math. Methods
Appl. Sci.}, \textbf{24} (2001), 405-426.

\bibitem{Smoller83} J. Smoller, Shock Waves and Reaction-Diffusion
Equations, Springer-Verlag, New York-Berlin, 1983.

\bibitem{Wang08} Z.A. Wang and T. Hillen, Shock formation in a chemotaxis model,
\textit{Math. Methods Appl. Sci.}, \textbf{31} (2008), 45-70.

\bibitem{Xin98} Z.P. Xin, Viscous boundary layers and their stability I., \textit{J. Partial
Differential Equations}, \textbf{11}(1998),  97-124.

\bibitem{Xin99} Z.P. Xin and T. Yanagisawa, Zero-viscosity limit of the linearized Navier-Stokes equations
 for a compressible viscous fluid in the
 half-plane, \textit{Comm. Pure Appl. Math.}, \textbf{52}(1999), 479-541.

\bibitem{Yang01} Y. Yang, H. Chen and W.A. Liu, On existence of global
solutions and blow-up to a system of reaction-diffusion equations
modelling chemotaxis, \textit{SIAM J. Math. Anal.}, \textbf{33}
(2001), 763-785.

\bibitem{Zhang07} M. Zhang and C.J. Zhu, Global existence of solutions to a
hyperbolic-parabolic system, \textit{Proc. Amer. Math. Soc.},
\textbf{135} (2007), 1017-1027.

\end{thebibliography}
\end{document}